\documentclass[english]{amsart}

\usepackage{babel}
\usepackage{amstext}
\usepackage{amsmath}
\usepackage{amsfonts}
\usepackage{amssymb}
\usepackage{latexsym}
\usepackage{ifthen}
\usepackage[all,2cell,dvips]{xy} \UseAllTwocells \SilentMatrices
\usepackage{xy}
\xyoption{all}
\renewcommand{\O}{{\mathcal O}}

\newcommand{\codim}{{\rm codim}}

\newtheorem{lemma1}{}[section]

\newenvironment{lemma}{\begin{lemma1}{\bf Lemma.}}{\end{lemma1}}
\newenvironment{example}{\begin{lemma1}{\bf Example.}\rm}{\end{lemma1}}
\newenvironment{question}{\begin{lemma1}{\bf Question.}\rm}{\end{lemma1}}

\newenvironment{theorem}{\begin{lemma1}{\bf Theorem.}}{\end{lemma1}}
\newenvironment{proposition}{\begin{lemma1}{\bf Proposition.}}{\end{lemma1}}
\newenvironment{corollary}{\begin{lemma1}{\bf Corollary.}}{\end{lemma1}}
\newenvironment{remark}{\begin{lemma1}{\bf Remark.}\rm}{\end{lemma1}}
\newenvironment{definition}{\begin{lemma1}{\bf Definition.}}{\end{lemma1}}

\newenvironment{conjecture}{\begin {lemma1}{\bf Conjecture.}}{\end{lemma1}}

\newcommand{\Q}{\ensuremath{\mathbb{Q}}}
\newcommand{\Z}{\ensuremath{\mathbb{Z}}}
\newcommand{\C}{\ensuremath{\mathbb{C}}}
\newcommand{\N}{\ensuremath{\mathbb{N}}}
\newcommand{\PP}{\ensuremath{\mathbb{P}}}
\newcommand{\D}{\ensuremath{\mathbb{D}}}

\newcommand{\merom}[3]{\ensuremath{#1:#2 \dashrightarrow #3}}

\newcommand{\holom}[3]{\ensuremath{#1:#2  \rightarrow #3}}

\makeatletter
\ifnum\@ptsize=0  \fi \ifnum\@ptsize=2  \fi 

\newcommand\sF{{\mathcal F}}
\newcommand\sG{{\mathcal G}}

\newcommand\sM{{\mathcal M}}
\newcommand\sO{{\mathcal O}}

\newcommand\sX{{\mathcal X}}

\newcommand\bZ{{\mathbb Z}}

\newcommand\bQ{{\mathbb Q}}

\DeclareMathOperator*{\Sym}{Sym}
\DeclareMathOperator*{\Tor}{Tor}

\newcommand\pure{\lq pure\rq \ }
\newcommand\canonical{\lq canonical\rq \ }
\newcommand\orbifold{\lq orbifold\rq \ }
\DeclareMathOperator*{\sat}{sat}
\DeclareMathOperator*{\orb}{orb}

\title {Special Orbifolds and Birational Classification: a Survey} 
\author{Fr{\'e}d{\'e}ric Campana}
\address{Fr\'ed\'eric Campana}
\email{campana@iecn.u-nancy.fr}
\date{\today}

\begin{document}

\begin{center}
{\Large }
\end{center}

\maketitle

\tableofcontents


\section{\bf Introduction: the Decomposition Problem.}

This text is an expansion of a talk given at the Schiermonnikoog Conference on Birational Classification, organised in May 2009 by C. Faber, E. Looijenga and G. van der Geer. It is intended to expose the main content of \cite{Ca04} and \cite{Ca07} (with some addtional topics in \S\ref{reldiff}, \S\ref{abund} or developments, as in \S\ref{fcpm})), as briefly and simply as possible, but essentially skipping the proofs. Some topics studied in \cite{Ca04} \and \cite{Ca07} have not been included here (such as orbifold versions of fundamental groups , universal covers, and function fields).

We shall show how to decompose, by functorial and canonical fibrations, arbitrary $n$-dimensional complex projective\footnote{Although the geometric results apply to compact K\" ahler manifolds without change, we consider here for simplicity this special case only.} varieties $X$ into
varieties (or rather \lq geometric orbifolds\rq $ )$ of one of the three \pure geometries determined by the `sign' (negative, zero, or positive) of the canonical bundle. These decompositions being birationally invariant, birational versions of these \pure geometries, based on the \canonical (or \lq Kodaira\rq $)$ dimension will be considered, rather.

This \pure trichotomy refines a more fundamental new dichtomy: `special' opposed to `general type' (the usual notion). `Special' turns out to be a suitable \orbifold combination of the first two pure geometries (canonical bundle negative or zero). More precisely, a variety $X$ is first decomposed\footnote{Everything works for orbifolds $(X\vert \Delta)$ as well.}, via a birationally unique fibration $c:X\to C(X)$ (its `core') into its antithetical parts: `special' (the fibres) and `general type' (the `orbifold base')\footnote{The term `orbifold' is used by Deligne-Mostow in almost the same sense as here. See footnote in \ref{dm} below. The term `orbifold pair' might possibly make the link to LMMP clearer.}. This is the new and most fundamental result of \cite{Ca04} and of the present text. 

Next, the `core' can be decomposed\footnote {Using an \orbifold version of Iitaka's $C_{n,m}$ conjecture.} as a composition $c=(M\circ r)^n$ of orbifold versions of the weak `rational quotient'\footnote{Also called `MRC fibration'.} $r$, and the Moishezon-Iitaka fibration $M$. When $X$ is \lq special\rq $, $ $C(X)$ is a point, and $X$ is thus a tower of fibrations with fibres of the first two geometries.

A crucial feature of these decompositions is indeed that, in order to deal with multiple fibres of fibrations, they need to take place in the larger category of `geometric orbifolds' $(X\vert \Delta)$. These are `virtual ramified covers' of varieties, which `virtually eliminate' multiple fibres of fibrations. Although formally the same as the `pairs' of the LMMP (see \cite{kmm}, \cite{KM}, \cite{BCHM} and the references there), they are here fully geometric objects equipped with the usual geometric invariants of varieties, such as sheaves of (symmetric) differential forms, fundamental group, Kobayashi pseudometric, integral points, morphisms and rational maps.

One expects theorems on birational properties of projective varieties to extend, together with their proofs (modulo some \orbifold adaptations) to the larger \orbifold category. This is illustrated by several examples below, the most important one being the weak positivity of direct images of pluricanonical sheaves in the \orbifold context. These \orbifold extensions are essential in the applications.

In the \orbifold category (but {\it not} in the category of varieties), most of the basic properties of the first two geometries are expected to be preserved by \lq extensions\rq $ $(ie: by the total space of a fibration, if satisfied by both the base and the general fibre). This naturally leads to conjecture, among many other things, that \lq special\rq orbifolds have an almost abelian fundamental group, and are {\it exactly} the ones with vanishing Kobayashi pseudometric or potentially dense set of rational points if defined over a number field. Also, conjecturally, the \lq core\rq $ $ should split the arithmetic and \C-hyperbolic properties of \lq orbifolds\rq $ $ in their two antithetical parts (\lq special\rq $ $ vs \lq general type\rq).

An implicit underlying theme of the considerations below is the `birational stability' of the cotangent bundle of a projective manifold. The question, vaguely formulated, is: to which extent is the amount of holomorphic sections of the  symmetric powers on $\Omega^p_X$ controlled by those of the pluricanonical bundles? Precise formulations are given in \S\ref{kvso}. THe expected answer is that these cotangent bundles are birationaaly stable, unless the the manifold is uniruled, in which case the birational stability is restored on the `rational quotient' of $X$. A similar answer is expected for orbifolds (see conjecture \ref{cok++}). We show in \S\ref{abund} how to deduce these conjectures from standard conjectures of the LMMP.

In this survey, many technicalities concerning the orbifold category have been skipped. The details and proofs are in \cite{Ca04} and \cite{Ca07}, where other related themes are treated. The results presented below on orbifold birational equivalence need to be completed on several essential points, and extended. See remarks \ref{rbir} and \ref{rbir'}. Several possibly accessible foundational questions remain unanswered. Essentially all of the definitions and results here apply to `smooth orbifolds' only. It is essential to extend these to the larger category of `log-canonical' pairs of the LMMP.

\

The investigations described below have been deeply influenced by K. Ueno's book [U75].

\section{\bf The three \lq pure\rq \ geometries}

We will now introduce the three \pure geometries which define the \lq elementary\rq $ $ objects 
of our classification theory. These \pure geometries are defined in two versions: first, `numerical', according to the `sign' of the canonical bundle, then, `birational', according to the refined `canonical' dimension. The definitions can be extended to \lq smooth geometric orbifolds\rq $ $ in the obvious way, but most known properties have not be shown to extend to this broader situation.

\subsection{The Trichotomy: numerical version}

The canonical divisor $K_X$ of a manifold $X$ has emerged as the major invariant of classification\footnote{Only recently it seems. First in the consideration of Enriques plurigenera and in the classification of surfaces by Kodaira and Shafarevich et al., and then in the three-dimensional classification initiated by Mori, Kawamata and Shokurov.}. Even to such an extent that one could almost define algebraic geometry as aiming at ``converting positivity (or negativity) of $K_X$ into geometry, topology, arithmetics of $X$". However, when the canonical bundle has generically negative directions, the canonical bundle is not sufficient, by itself, to control all directions of the cotangent bundle. See \S \ref{kvso} below.

The three elementary classes are thus the cases where $K_X$ is either antiample, or (numerically)
trivial, or ample which we will respectively denote by: $K_X<0$, $K_X \equiv 0$ and $K_X >0$.

For these three classes of manifolds some of the qualitative properties of curves are known to hold in higher dimension, the others being conjectural, which is indicated by the sign ``?". The symbol ``??" means that no general conjecture is presently even formulated, while `essentially' means: `outside of a proper algebraic subset'. We denote with $X^*(k)$ the Zariski closure of the set of $k$-rational points of $X$.

\[
\begin{tabular}{|l|l|l|l|}
\hline

$K_X$ & $\pi_1(X)$ & $d_X$ & $X^*(k)$
\\
\hline
$<0$ & $\{ 1 \}$ & $\equiv 0$ & $X ?$
\\
\hline
$\equiv 0$ & almost abelian & $\equiv 0 ?$ & $X ?$
\\
\hline
\hline
$>0$ & ??& `essentially' a metric ? & `essentially' finite ?
\\
\hline
\end{tabular}
\]

We see that essentially nothing is known in higher dimensions on the arithmetic side, and also for varieties of general type. The vanishing of the Kobayashi pseudometric on a Fano manifold follows from their rational connectedness (see below), which is one of the main applications of Mori theory (\cite{Ca92}, \cite{KMM92}).

\subsection{The Trichotomy: birational version}

We define now the three pure geometries as follows (except for 1, this is Ueno's trichotomy in \cite{U75}, \S 11):

1. $\kappa_+(X)=-\infty$ (see definition \ref{dk+} below), which is conjecturally equivalent (see conjecture \ref{ck}) to rational connectedness.

2.  $\kappa(X)=0$,

3.  $\kappa(X)=\dim X$.

In analogy to the case of curves and the numerical version of the trichotomy, 
the expected properties (all of them conjectural)
are as follows:

\[
\begin{tabular}{|l|l|l|l|}
\hline
$\kappa$ & $\pi_1(X)$ & $d_X$ & $X^*(k)$
\\
\hline
$\kappa_+(X)=-\infty$ & $\{ 1 \}?$ & $\equiv 0 ?$ & $X ?$
\\
\hline
$\kappa(X)=0$ & almost abelian? & $\equiv 0 ?$ & $X ?$
\\
\hline
$\kappa(X)=\dim X$ & & essentially $>0$ ? & essentially finite ?
\\
\hline
\end{tabular}
\]

\

\begin{remark} Let us give some reasons supporting these conjectures.

1. If the conjecture \ref{ck+} holds, a manifold with $\kappa_+(X)=-\infty$
is rationally connected, so its Kobayashi pseudometric vanishes, and $X$ is simply connected.

2. If $\kappa(X)=0$, the (conjectural) minimal model program
implies that $X$ is birational to a mildly singular\footnote{Technically: `terminal'.} variety $X'$ such that
the canonical divisor $K_{X'}$ is numerically trivial. Singular versions
of the Bogomolov decomposition theorem should imply that the fundamental group is 
almost abelian. The ``classification" of varieties with numerically trivial canonical
divisor suggests that the Kobayashi pseudometric vanishes.  For curves and surfaces this is known. In the case of ``Hyperk\" ahler" manifolds, only weaker versions are known, using twistor spaces. But the proof suggests this vanishing. The case of Calabi-Yau manifolds (simply-connected, with $h^{(0,2)}=0)$ is much more conjectural. The computation of Gromov-Witten invariants on some examples shows however the existence of infinitely many rational curves.
\end{remark}

\subsection{A remark on the LMMP}

\

The aim of the MMP is actually to reduce the `birational' version above to the `numerical' version by constructing `minimal models', and converting hypothesis of pseudo-effectivity for adjoint line bundles by nefness first, and then by semi-ampleness. This task naturally leads (for reasons apparently different from the ones here)\footnote{Although the proofs of Kawamata-Viehweg Theorem, as well as their approach to weak positivity of direct images of pluricanonical sheaves rely on cyclic covers, the same being implicitely true of  Kawamata's subajunction theorem.} to introduce the very same `pairs' $(X\vert \Delta)$ as below, which permits an inductive treatment on the dimension, by producing `log-canonical centers', and to extend the MMP to the LMMP. However, while the LMMP considers only the canonical bundle $K_X+\Delta$ of pairs, our approach leads us to equip them naturally with many other geometric invariants.

\subsection{Rational connectedness and $\kappa_+=-\infty$.}

\

\

The following is one of the central problems in birational classification, and the first step of the so-called `Abundance conjecture':

\begin{conjecture}
\label{ck}
Let $X$ be a projective manifold such that $\kappa(X)=-\infty$. Then $X$ is uniruled.
\end{conjecture}

This conjecture is known up to dimension $3$ by the already classical work of S. Mori, Y. Kawamata, S. Shokurov and Y. Miyaoka.

The class of projective manifolds such that $\kappa(X)=-\infty$ contains all products $\PP^1 \times Y$, and so does not define
a \lq pure\rq$ $ geometry . Hence the following definition:

\begin{definition}\label{dk+} We define: $
[\kappa_+(X)=-\infty]  :\Leftrightarrow [\kappa(Y)=-\infty,$ for all positive dimensional manifolds $Y$ dominated by $X]$,
where: ``$Y$ is dominated by $X"$ means that there exists a dominant rational map \merom{g}{X}{Y}.

More generally, for any, $X$, we define $\kappa_+(X):=max\{\kappa(Y)\}$, for all positive-dimensional $Y$ dominated by $X$.
\end{definition}

Obviously, $n\geq\kappa_+(X)\geq \kappa(X)$, for any $X$, and $\kappa_+(X)=-\infty$ if $X$ is rationally connected. Conversely:

\begin{conjecture}
\label{ck+}
A manifold $X$ is rationally connected if 
$
\kappa_+(X)=-\infty.
$
\end{conjecture}

Recall the existence of the following map $r$, the
`rational quotient' (or `MRC fibration' in \cite{KMM92}), which describes the $\kappa_+=-\infty$ `part' of any $X$.

\begin{theorem}\label{rq} \cite{Ca81, KMM92, GHS}
Let $X$ be a projective manifold, then there exists a unique rational map:
\merom{r_X}{X}{R(X)} such that: 
\begin{itemize}
\item its fibres are rationally connected,
\item R(X) is not uniruled.
\end{itemize}
\end{theorem}

{\bf Remark.} It easily follows from this theorem that conjecture \ref{ck}
implies conjecture \ref{ck+}. Indeed, if $r_X:X\dasharrow R(X)$ is not the constant map, $R(X)$ is not uniruled, thus $\kappa(R(X))\geq 0$ by conjecture \ref{ck}, contradicting the assumption $\kappa_+=-\infty$ $\square$

\subsection{Birational stability of  Cotangent bundles}\label{kvso}

\

\

Let $T(X)$ be the (birational) complex algebra of all `covariant' holomorphic tensors (ie: sections of $\otimes \Omega_X^1)$ on the complex projective manifold $X$. This algebra contains fundamental information on the rational fibrations $f:X\dasharrow Y$ (we give two examples below). The study of this algebra is difficult, since it is deduced from the rank-$n$ cotangent bundle of $X$. It is thus of central importance to deduce its qualitative structure from positivity properties of the rank-one canonical bundle alone, or in other words, to establish the `birational stability' of $\Omega^1_X$ (in a precise sense given below). This is true in the numerical version, i.e when the numerical properties of $K_X$ are considered (see below). This fails however completely when the birational invariants are considered, specifically when $\kappa(X)=-\infty$. This property is much too weak to control all directions of the (co)tangent bundle, as shown by $X=\Bbb P^1\times Y, Y$ any projective manifold. This stability should hold, however, under the refined condition $\kappa_+(X)=-\infty$, which considers the canonical bundles of all `quotients' of $X$, such as $Y$ when $X=\Bbb P^1\times Y$.

When $X$ is of `pure geometry' in the `numerical version', the canonical bundle determines the structure of $T(X)$. Indeed, when $X$ is Fano, and so rationally connected, this algebra reduces to $\Bbb C$, the constants. When $K_X\equiv 0$, Miyaoka generic semi-positivity\footnote{Or, alternatively, using Ricci-flat metrics, the parallelism of these tensors.} implies that non-zero elements of $T(X)$ do not vanish anywhere, which permits the explicit determination of $T(X)$ by representation theory, and implies that this algebra is finitely generated (See \cite{Ca95}, or \cite{Pe} for the details of the argument). When $K_X$ is ample, the situation is not uniform, the symmetric part of the algebra above being possibly reduced to the constants (for hypersurfaces of the projective space, when $n\geq 2)$. But sections of higher jets bundles are expected to always exist in abundance; this is an essential issue in hyperbolicity questions. See, for example \cite{D},\cite{DMR}, where many more references are also given.

The proofs of these results, which permit to control the rank-$n$ bundle $\Omega_X^1$ by its determinant line bundle, use deep and indirect results: either Mori's reduction to char $p>0$, or Yau's solution of the Calabi Conjecture.

For general $X$ of intermediate positivity, these algebras (which are birational invariants, possibly even also deformation invariants) have apparently not been investigated.

For the `pure geometries' in their `birational versions', essentially nothing is proved, although the same behaviour as in the `numerical' version is expected. Indeed, if conjecture \ref{ck+} holds, then $T(X)=\Bbb C$ when $\kappa_+(X)=-\infty$, as well. When $\kappa(X)=0$, the Abundance conjecture implies as above that $T(X)$ should have the same structure as when $K_X\equiv 0$ (see \cite{Ca95}). 

Bounding the positivity of $\Omega_X^1$ also permits to control the fibrations on $X$ and other invariants, such even as $\pi_1(X)$. More precisely, let us formulate in the conjectures \ref{ck++} and \ref{ck^+} below the expected `birational stability' of $\Omega^1_X$:

\begin{definition}\label{dk++} Let $X$ be a projective complex connected manifold. Then: $\kappa_{++}(X):=max_{\{L,p>0\}}\{\kappa(X,L)\}$, $L\subset \Omega_X^p$ being a  rank $1$ coherent subsheaf $\}$.
\end{definition}

Obviously: $n\geq \kappa_{++}(X)\geq\kappa_+(X)\geq \kappa(X)$.

\begin{conjecture} \label{ck++} If $\kappa(X)\geq 0$, then $\kappa_{++}(X)=\kappa(X)$.

If $\kappa_+(X)=-\infty$, then $\kappa_{++}(X)=-\infty$.

More precisely: $\kappa_{++}(X)=\kappa(R(X))$, for any $X$, $R(X)$ being its `rational quotient'. In particular: $\kappa_{++}(X)=\kappa_+(X)$, for any $X$.
\end{conjecture}

\begin{remark}\label{ck0'} Easy arguments (using the Moishezon-Iitaka fibration) show that if the conjecture holds for $\kappa=0$, it also holds when $\kappa\geq 0$. The same arguments also apply in the orbifold case (see the conjecture \ref{cok++}).
\end{remark}

 If $f:X\dasharrow Y$ is a nonconstant fibration, $\kappa(Y)=\kappa(X,f^*(K_Y))\leq \kappa_{++}(X)$.

The above conjecture implies thus: $\kappa(Y)\leq 0$ if $\kappa(X)=0$ (which follows from Iitaka's $C_{n,m}$ too).

Bogomolov's theorem (\cite{Bo78}) asserts that $\kappa(X,L)\leq p$ if $L\subset \Omega_X^p$ is a  rank $1$ coherent subsheaf, and that $L=f^*(K_Y)$ generically on $X$, for some unique $f:X\dasharrow Y$, if equality holds. However, $Y$ needs not then be of general type. The difference between $L$ and $f^*(K_Y)$ lies in the multiple fibres of $f$, and the `orbifold base' $(Y\vert \Delta_f)$ is of general type in this case, as will be seen below (in the more general `orbifold' context).

The very same results and expectations hold also for `smooth geometric orbifolds' (see the conjecture \ref{cok++} below).

\

$\bullet$ Another similar birational invariant which has implications on the structure of the universal cover is:

\begin{definition}\label{dk^+}\cite{Ca95}  Let $X$ be a projective complex connected manifold. Then: $\kappa^+(X):=max_{\{\sF,p>0\}}\{\kappa(X,det(\sF))\}$, where $\sF\subset \Omega_X^p$ is a coherent subsheaf.
\end{definition}

Similarly to \ref{ck++}:

\begin{conjecture}\cite{CA95} \label{ck^+} If $\kappa(X)\geq 0$, then $\kappa^{+}(X)=\kappa(X)$.

If $\kappa^+(X)=-\infty$, then $X$ is rationally connected.

More precisely: $\kappa^+(X)=\kappa(R(X))$, for any $X$.
\end{conjecture}

This conjecture is established wenn $\kappa^+(X)=n$ (see \cite{CP05}). When $\kappa(X)=0$, it implies Ueno's Conjecture that $h^0(X,\Omega_X^p)\leq (_p^n)$ if $\kappa(X)=0$, which is still open.

Using $L^2$-theory, Atiyah's index theorem, and Gromov's Poincar\'e series \cite{Gr92}, it is shown in \cite{Ca95} that:

\begin{theorem} Assume that $\chi(\sO$$_X)\neq 0$. Then: $\gamma(X)\leq \kappa^+(X)$.
\end{theorem}

Here $\gamma(X):=dim(\Gamma(X)$, where $\gamma_X:X\dasharrow \Gamma(X)$ is the `Gamma-reduction' (or `Shafarevich map') of $X$, with general fibres the largest subvarieties in $X$ whose fundamental group has finite image in $\pi_1(X)$. The fibre-dimension of $\gamma_X$ is also the dimension of the largest compact connected analytic subset through a general point of the universal cover of $X$.

This result implies in particular that: $\pi(X)=\{1\}$ if $\kappa^+(X)=-\infty$, that $\pi_1(X)$ is finite if $\kappa^+(X)=0$ and $\chi(\sO$$_X)\neq 0$, and that $X$ is of general type if $\gamma(X)=n$ and $\chi(\sO$$_X)\neq 0$ (\cite{CP05}). Observe that Abelian varieties have $\kappa=\kappa^+=n$, and $\gamma(X)=n$, so that the condition $\chi(\sO$$_X)\neq 0$ is essential in the last two statements. An example of a threefold $X$ of general type with $\gamma(X)=n=3$, and $\chi(\sO$$_X)=0$ is given in \cite{EL}.

\subsection{ The decomposition problem: a failed attempt.}\label{dpfa}

\

\

Any curve belongs to one of the three pure geometries.
In dimension at least two, this is no longer true, and fibrations (ie: dominant rational maps  \merom{f}{X}{Y}
with connected general fibre)
are needed to decompose arbitrary manifolds $X$
into pieces of `pure geometry'. 
\

The two classical maps: first $r:X\to R(X)$ (the ``rational quotient" or ``MRC"-fibration, see ), and $M=M_X:X\to M(X)$ (the Moishezon-Iitaka fibration, if $\kappa(X)\geq 0)$ seem, at first sight, to provide a solution to this decomposition problem. Indeed, the first one eliminates the $RC$ `part' of $X$, and the second its $\kappa=0$ `part'. 

Since $r$ has a non-uniruled base $R(X)$, it has conjecturally $\kappa(R(X))\geq 0$ (by conjecture \ref{ck}). So that $M\circ r:X\to M(R(X))$ should be well-defined for any $X$. However, $M(R(X))$ is not of general type in general, and one needs to iterate and consider $(M\circ r)^n:X\to MR^n(X)$ to reach a base $MR^n(X)$ of general type (possibly a point). The problem is that `parts' of general type can be hidden in the seemingly `general type-free' fibres of $(M\circ r)^n$.The simplest example is given in \ref{exampleiitaka} below. This $(M\circ r)$-process thus reduces the dimension, but does not describe the structure of general $X's$. 

Notice that the LMMP, which aims at producing a `numerical' version of the maps $M$ and $r$ (converting, for appropriate adjoint $\Bbb Q$-line bundles, the word `pseudoeffective' into `nef' first , and next into `semi-ample') does not consider this question.

As we shall see, this decomposition process works, but only in the `orbifold category' defined below.  Indeed, as shown by the example \ref{exampleiitaka} below, the main problem is that the couple $(X_y,Y)$, where $X_y$ is the general
fibre and $Y$ the base of some fibration $f:X\to Y$, {\it do not determine} the qualitative geometry of $X$, even when $f$ is one of
the natural fibrations of the minimal model (or of the Moishezon-Iitaka) program. This failure is due to the presence of multiple fibres, and disappears in the category of `geometric orbifolds', which consists precisely in encoding this data. Surprisingly, this single addition is sufficient to correct the above $(M\circ r)^n$ decomposition. Working in the category of geometric orbifolds is thus, not only necessary, but also sufficient to solve the decomposition problem.

Let us illustrate the problem at hand with the following very simple example \ref{exampleiitaka}.

\begin{example}
\label{exampleiitaka}
Let $C$ be a hyperelliptic curve, and let $E$ be an elliptic curve (both defined over some number field $k$, say).
Let \holom{i}{C}{C} be the hyperelliptic involution that exchanges the fibres of the morphism $C \rightarrow \PP^1$ 
and let $\holom{t}{E}{E}$ be a translation of $E$ of order $2$, thus $t$ has no fixed point. 

The quotient map
$u:X':=E\times C\to 
X := (C \times E)/< i \times t>
$ is thus \'etale, and  $\kappa(X)=1$.

We get so a commutative
diagram, in which the vertical arrows are also the Moishezon-Iitaka fibrations:
\[
\xymatrix @C=3.5cm { X'=C \times E \ar[d]^{M_{X'}} \ar[r]^{u}_{(2:1) \ \mbox{\footnotesize \'etale}} & X= X'/<i\times t>
\ar[d]^{M_X}
\\
C \ar[r]^{v}_{(2:1) \ \mbox{\footnotesize ramified}} & \PP^1=C/<i> }
\]
\end{example}

Because $\pi_1(X)$, $d_X$ and $X(k)$ are essentially invariant under finite unramified covers, these invariants coincide essentially with those of $X'$, and thus radically differ from those expected by considering the couple $(X_y,Y)=(E,\PP^1)$ of generic fibre and base of $M_X:X\to \PP^1$. Seeing this fibration as a `twisted product' $E \times \PP^1$, 
one would expect $X$ to have
an almost abelian fundamental group, $d_X\equiv 0$, and $X(k)$ Zariki dense (after a finite extension).

The qualitative geometry of $X$ can be recovered by looking at $X$ only, since $X'$ is nothing, but the normalisation of the fibre product $X\times _{\Bbb P^1}C$, and $C$ is simply a ramified cover of $\Bbb P^1$ which ramifies at order two exactly over the points $p\in \Bbb P^1$ over which $M_X$ has a double fibre. Under the above normalised fibre product, these multiple fibres are thus eliminated.

We will now generalise this construction. Let $f:X\to Y$ be a fibration. Let $\Delta_f$ be the $\bQ$-divisor on $Y$ uniquely defined by the multiple fibres of $f$, in such a way that base-changing $f$ by a {\it local} finite cover $Y'\to Y$ ramifying `exactly' over $\Delta_f$ eliminates in codimension one the  multiple fibres of the fibration $f:X\to Y$. Since a finite ramified cover $Y'\to Y$ such as $C\to \Bbb P^1$ above  does not exist in general globally, we need to work directly with the pair $(Y\vert \Delta_f)$.

\section{\bf Geometric Orbifolds}

We give only the definitions needed for the considerations below. The other definitions are in \cite{Ca07}.

Throughout this section we will denote by $X$ a smooth projective manifold of dimension $n$.

\subsection{Geometric Orbifolds}

\begin{definition}\label{dm}
Let $Y$ be a normal connected complex projective variety. An orbifold divisor is an effective $\Q$-divisor 
$\Delta = \sum_{D \subset Y} (1- \frac{1}{m(D)}) D,$
where: 
\begin{itemize}
\item the sum ranges over all prime divisors on $Y$, 
\item $m(D) \in (\Q\cap [1,+\infty[) \cup \{ +\infty \}$,
\item $m(D)=1$ for all but a finite number of $D$'s.
\end{itemize}
A `geometric orbifold' (or simply: `orbifold')\footnote{In \cite[\S14, pp. 135-141]{DM}, this term is employed in a sense which is essentially equivalent to ours in the finite, integral, smooth case. This reference was pointed to me by F. Catanese.} is a couple $(Y\vert \Delta)$ where $Y$ is smooth\footnote{Being normal and $\bQ$-factorial is actually sufficient for most of the definitions given here.} projective and $\Delta$ an orbifold divisor on $Y$; it is said to be finite (resp. integral, resp. logarithmic) if for all $D\subset Y$, one has: $m(D)<+\infty$ (resp. $m(D)\in (\bZ\cup +\infty)$, resp. $m(D)=+\infty)$.

The `support' $Supp(\Delta)=\lceil\Delta\rceil$ of $\Delta$ is the (finite) union of all $D's$ such that $m(D)>1$.

The geometric orbifold $(Y\vert \Delta)$ is said to be `smooth' if $Y$ is smooth, and if $Supp(\Delta)$ is a divisor of normal crossings.

When $\Delta=0$ (resp. when $\Delta$ is logarithmic) the orbifold $(X\vert \Delta)$ is identified with $X$ (resp. with the quasi-projective variety $U:=X-\Delta)$.

We define a lattice order on the set of orbifold divisors on $X$ by writing $\Delta'\geq \Delta$ if $(\Delta'-\Delta)$ is effective.

\end{definition}

Geometric orbifolds $(Y\vert \Delta)$ with finite multiplicities interpolate between proper or compact orbifolds (when $\Delta=0)$, and open orbifolds (when the multiplicities are all infinite). When moreover integral, they may be considered as virtual covers of $Y$ ramified over each divisor $D$ of $Y$ with multiplicity $m_{\Delta}(D):=m(D)$ in the notation above.

Geometric orbifolds thus coincide with the `log-pairs' (with rational coefficients) of the LMMP. Notice that a smooth orbifold is log-canonical, and klt if finite.

The origin and main source of geometric orbifolds here are the `orbifold bases' of fibrations (see \S \ref{ob} below).

\

\subsection{Orbifold Invariants}\label{inva}

\

\

The most fundamental one is the following:

\begin{definition} Assume $Y$ to be smooth (or $K_Y$ to be $\Q$-Cartier).
The canonical bundle of the orbifold $(Y\vert \Delta)$ is defined as
$K_{Y\vert\Delta}=K_Y + \Delta$,
and the canonical dimension of $(Y,\Delta)$ is 
$
\kappa(Y\vert\Delta):=\kappa(Y, K_{Y\vert\Delta}).$
\end{definition}

{\bf Remark.} One thus has: $\dim Y \geq \kappa(Y\vert\Delta) \geq \kappa(Y).$

\

Other important invariants are, when $(X\vert \Delta)$ is smooth:

\

$\bullet$ The locally free sheaves $S^N\Omega^q(Y\vert\Delta)$ for any $N,q\geq 0$, when $(Y\vert\Delta)$ is smooth, and $p$-dimensional. They agree with $Sym^N(\Omega^q_Y)$ if $\Delta=0$, and with $Sym^N(\Omega^q_Y(log (D)))$ if $\Delta=Supp(\Delta)=D$ is logarithmic. In general they interpolate between these two cases. 

In local analytic `adapted' coordinates $z=(z_1,\dots,z_p)$, in which the support of $\Delta$ consists of coordinates hyperplanes, the sheaf $S^N\Omega^q(Y\vert \Delta)$ is generated as an $\sO_Y$-module, by the elements $\frac{dz_{(J)}}{z_{(J)}}:=z^{\lceil k/m\rceil}.\otimes _{\ell=1}^{\ell=N}\frac{dz_{(J_{\ell})}}{z_{(J_{\ell})}}$, where: 

1. the $J_{\ell}$ are increasing subsets of $\{1,\dots,p\}$ of cardinality $q$, lexicographically ordered in incresasing order, and $\frac{dz_{(J_{\ell})}}{z_{(J_{\ell})}}:=\wedge_{j\in J_{\ell}}\frac{dz_j}{z_j}$.

2. $k_j$ is the number of $\ell$'s such that $j\in J_{\ell}$, for each $j\in \{1,\dots,p\}$.

3. $z^{\lceil k/m\rceil}:=\Pi_{j=1}^{j=p}z_j^{\lceil k_j/m_j\rceil}$, where $m_j$ is the $\Delta$-multiplicity of the hyperplane of equation $z_j=0$.

Note that $\frac{dz_{(J)}}{z_{(J)}}:=\Pi_{j=1}^{j=p}(z_j^{-\lfloor k_j.(1-\frac{1}{m_j})\rfloor}).\otimes _{\ell=1}^{\ell=N}dz_{(J_{\ell})}$, for any $(J)$.

\

The following lemma will be used in the proof of  lemma \ref{fanok++}.

\begin{lemma}\label{diffsheaf} With the notations above, assume that, for any $j\in \{1,\dots,p\}$, one has: $m_j>1$. Then, for any index $(J)=(J_1,\dots,J_q)$ above, there exists some $j=j((J))$ such that  $\lfloor k_j.(1-\frac{1}{m_j})\rfloor>0$ if $N\geq \frac{p}{q.(1-1/m)}$, where $m:=inf_j\{m_j\}>1$.
\end{lemma}

{\bf Proof:} Assume, by contradiction, that $k_j.(1-1/m)<1$ for each $j$. Since $\sum_jk_j=Nq$, we get: $Nq.(1-1/m)\leq \sum_jk_j.(1-1/m_j)<p$, and the claim $\square$

\

$\bullet$ The sheaves of holomorphic tensors of type $T^r_s$ on $(Y\vert \Delta)$ are defined similarly. See \cite{Ca07} for more details, and \cite{JK} for some additional properties.

\

In the proof of lemma \ref{desc}, we shall need a relative version, of independent interest, of the sheaves of differentials defined above. See \S\ref{reldiff} below.

\

 When $\Delta$ is integral, one can also naturally define (see \cite{Ca07} for details):

 $\bullet$ The fundamental group: $\pi_1(Y\vert\Delta)$, and the universal cover when $(Y\vert \Delta)$ is moreover smooth.
 
 $\bullet$ The Kobayashi pseudometric $d_{Y\vert\Delta}$ (see example \ref{d} below)
 
 $\bullet$ The notion of set of integral points $(Y\vert\Delta)(\O$$_{k,S})$ if $k$ is a number field over which $(Y\vert\Delta)$ is defined, with ring of integers $\O$$_{k}$, and a finite set of places $S$. (This depends on the choice of a model of $(Y\vert\Delta)$ over $\O$$_{k})$. See \cite{Ab} for a detailed exposition.

\subsection{Orbifold Morphisms.}

\

\

Let $f:X\to Y$ be a regular map between projective normal varieties. Assume $Y$ to be smooth (or $\bQ$-Cartier). Let $\Delta_X$ and $\Delta_Y$ be orbifold divisors on $X$ and $Y$ respectively.

For any prime Weil divisor $E$ on $Y$, let: $f^*(E):=\sum_{D\subset X} t_{E,D}.D$, the sum running over all prime divisors on $X$. Thus $t_{E,D}>0$ if and only if $f(D)\subset E$. 

We write $m_X(D):=m_{\Delta_X}(D)$, and similarly for $E,Y$.

\begin{definition} We say that $f:(X\vert\Delta_X)\to (Y\vert\Delta_Y)$ is an orbifold morphism if:

\begin{itemize}
\item $f(X)$ is not contained in $Supp(\Delta_Y)$.
\item For any $E,D$ as above,  if $t_{E,D}>0$, then: $t_{E,D}.m_X(D)\geq m_Y(E)$.
\end{itemize}

We shall say that the orbifold morphism $f:(X\vert\Delta_X)\to (Y\vert\Delta_Y)$ is a `classical' orbifold morphism if $(Y\vert\Delta_Y)$ and $(X\vert\Delta_X)$ are integral, and if moreover, in the second condition above, $``t_{E,D}.m_X(D)\geq m_Y(E)"$ is replaced by: $``m_Y(E)$ divides $t_{E,D}.m_X(D)"$:
\end{definition}

\begin{example}\label{d} Assume $(Y\vert \Delta_Y)$ is smooth. A holomorphic map $h:\D\to Y$ such that $f(\D)$ is not contained in $Supp(\Delta_Y)$ defines an orbifold morphism (resp. a classical orbifold morphism) $h:\D\to (Y\vert \Delta_Y)$ if, for any $z\in \D$ such that $f(z)\in Supp(\Delta_Y)$, the order of contact of $f(\D)$ with every branch $E$ of $Supp(\Delta_Y)$ at $f(z)$ is at least (resp. is a multiple of) $m_Y(E)$.(The notion of orbifold morphism does not need the projectivity of $X$ or $Y)$.

This leads to a definition of the Kobayashi pseudo-metric $d_{Y\vert \Delta_Y}$ for a smooth orbifold: this is the greatest pseudo-metric $\delta$ on $Y$ such that $h^*(\delta)\leq d_{\D}$ for any orbifold morphism $h:\D\to (Y\vert \Delta_Y)$, with $d_{\D}$ the Poincar\'e metric on $\D$. If $(Y\vert \Delta_Y)$ is integral, we can define similarly the `classical' Kobayashi pseudo-metric $d^*_{(Y\vert \Delta_Y)}$ on $(Y\vert \Delta_Y)$, by replacing the orbifold morphisms $h:\D\to (Y\vert \Delta_Y)$ above by their `classical' analogues. We have natural inequalities: $d^*_{(Y\vert \Delta_Y)}\geq d_{(Y\vert \Delta_Y)}\geq d_Y$.

When $Y$ is a curve, we have: $d^*_{(Y\vert \Delta_Y)}=d_{(Y\vert \Delta_Y)}$ (\cite{Rou06}, answering a question in \cite{CW}). 
\end{example}

The notions of morphisms and differential forms are compatible\footnote{When one considers orbifolds, even log-terminal, with $X$ being $\Bbb Q$-factorial, but not factorial, the definition of an orbifold morphism might need to take into account the non-factoriality of $X$, by introducing a suitable factor of local non-factoriality. The simplest example, pointed out by the referee, being the cone over the conic equipped with two lines of the ruling with multiplicity $2$ as $\Delta$, and blown up at the vertex, exhibits a difference between the usual orbifold theory, where the exceptional divisor has multiplicity $1$, and the definition above, which gives multiplicity $2$.}:

\begin{theorem}\label{mor}(\cite{Ca07}, 2.7 and 2.12) Let $f:X\to Y$ be holomorphic, with $X,Y$ smooth, and $\Delta_X,\Delta_Y$ as above. Assume that $f(X)$ is not contained in $Supp(\Delta_Y)$.

The following conditions are then equivalent:
\begin{itemize}
\item $f:(X\vert \Delta_X)\to (Y\vert \Delta_Y)$ is an orbifold morphism
\item $f^*(S^N\Omega^p_{(Y\vert \Delta_Y)})\subset S^N\Omega^p_{(X\vert \Delta_X)}, \forall p,N\geq 0$
\item If $(X\vert \Delta_X)$ is integral: for any orbifold morphism $h:\D\to (X\vert \Delta_X)$, $f\circ h:\D\to (Y\vert \Delta_Y)$ is also an orbifold morphism.

\end{itemize}
\end{theorem}

\

\subsection{Orbifold Birational maps.}\label{bir}

\

\

The following definitions need to be generalised to the larger class of log-canonical and klt orbifolds (instead of smooth ones).

\begin{definition} Let $f:X\to Y$ be a birational map between smooth projective manifolds equipped with orbifold divisors $\Delta_X,\Delta_Y$. We say that $f:(X\vert \Delta_X)\to (Y\vert \Delta_Y)$ is an elementary birational orbifold morphism (an `ebom') if:
\begin{itemize}
\item $(Y\vert \Delta_Y)$ and $(X\vert \Delta_X)$ are smooth
\item $f:(X\vert \Delta_X)\to (Y\vert \Delta_Y)$ is an orbifold morphism
\item $f_*(\Delta_X)=\Delta_Y$
\end{itemize}

The smooth orbifolds $(X\vert \Delta_X)$ and $(Y\vert \Delta_Y)$ are said to be birationally equivalent if they can be connected by a chain of `eboms'
\end{definition}

From theorem \ref{mor}, we deduce that the spaces of global sections of the sheaves $S^N(\Omega^p_{*})$ are orbifold birational invariants, since these sheaves are locally free.

\begin{remark}\label{rbir'} In general, two birational smooth orbifolds are, however, not birationally dominated by a third one (see \cite{Ca07}, 2.32 and 2.33). This is the source of considerable technical, even possibly conceptual, difficulties.
\end{remark}

\begin{remark} Two logarithmic orbifolds are birationally equivalent in the above sense if and only if their open parts $(X-\Delta)$ are (algebraically) isomorphic. Thus birational invariants in the orbifold category produce (many new) invariants of quasi-projective varieties.
\end{remark}

\subsection{Relative differentials.}\label{reldiff}

\

\

Let $f:X\to Y$ be a fibration (ie: a regular surjective and connected map) between complex projective connected manifolds. Let $\Delta$ be an orbifold divisor on $X$, such that $(X\vert \Delta)$ is smooth. Let $\Delta_Y:=\Delta_{f\vert \Delta}$ the associated base orbifold divisor on $Y$. By suitable blow-ups of $X$ and $Y$, and orbifold modifications of $(X\vert \Delta)$, we can assume that: 

1. $f$ is neat.

2. $(Y\vert \Delta_Y)$ is smooth.

3. $f:(X\vert \Delta)\to (Y\vert \Delta_Y)$ is an orbifold morphism.

\

From the property (3) above, we deduce, for any pair of nonegative integers, $N,q$, a natural map of sheaves $f^*:S^N(\Omega^q(Y\vert \Delta_Y))\to S^N(\Omega^q(X\vert \Delta))$.

We shall derive sufficient conditions for the above map to be surjective at the level of global sections. The conditions bear on the sheaves $S^N(\Omega^q(X_y\vert \Delta_y))$, where $X_y$ is a generic smooth fibre of $f$, while $\Delta_y$ is the restriction of $\Delta$ to $X_y$. Notice that Sard's (or Bertini's) theorem implies that $(X_y\vert \Delta_y)$ is a smooth orbifold.

\begin{proposition}\label{f*surj} In the above situation, $f^*:H^0(Y,S^N(\Omega^q(Y\vert \Delta_Y)))\to H^0(X,S^N(\Omega^q(X\vert \Delta)))$ is surjective if, for any finite sequence $(N_h,q_h),h=1,\dots, t$ of pairs of positive integers, $H^0(X_y,S^{N_1}\Omega^{q_1}(X_y\vert \Delta_y)\otimes \dots \otimes S^{N_t}\Omega^{q_t}(X_y\vert \Delta_y))=\{0\}$.
\end{proposition}

{\bf Proof:} $\bullet$ The first (and longer) step consists in studying the natural filtration of the sheaves of relative orbifold differentials. 

Let $X_y$ be a generic fibre as above, and let $M$ be its (trivial) normal bundle, of rank $(n-p)$. There is natural increasing filtration $F^s, s=-1,0,\dots, q$ of $\Omega^q(X)_{\vert X_y}$ with graded pieces equal to $F^s/F^{s-1}=\Omega^s(X_y)\otimes \Lambda^{q-s}(M)$, for $s=0,1,\dots, q$. 

This filtration induces an increasing filtration $G^tF^s$ of $Sym^N(F^s)$, for $t= -1,0,\dots,N$, associated to the exact sequence $0\to F^{s-1}\to F^s\to (F^s/F^{s-1})\to 0$. 

It has thus graded pieces $G^tF^s/G^{t-1}F^s=Sym^t (F^{s-1})\otimes Sym^{N-t}(F^s/F^{s-1})$.

Altogether we get a filtration of $Sym^N(\Omega^q(X))_{\vert X_y}$ with graded pieces $\otimes_{s=0}^{s=q}Sym^{N_s}(F^{s}/F^{s-1})=\otimes_{s=0}^{s=q} Sym^{N_s}(\Omega^{s}(X_y)\otimes \Lambda^{q-s+1}(M))$, for all multi-indices $(N_0,\dots,N_q)$ of sum $N$.

In local adapted coordinates $(z)=(z_1,\dots,z_n)$, with $\Delta$ supported on the coordinate hyperplanes $z_j, j=1,\dots,(n-p)$, and $f(z)=(y_1:=z_n-p+1,\dots,y_p:=z_n)$, recall that a basis of the $\sO_X$-module $S^N(\Omega^q(X\vert \Delta))$ is given by the following: 
$$ \frac{dz_{(J)}}{z_{(J)}}:=\frac{1}{z_{(J)}}.\otimes _{\ell=1}^{\ell=N}dz_{(J_{\ell})},$$ for $(J)=(J_1,\dots,J_N)$ any (lexicographically oredered) $N$-tuple of subsets $J_{\ell}$ of $\{1,2,\dots,n\}$ of cardinality $q$, where: $$z_{(J)}:=\Pi_{j=1}^{j=(n-p)}(z_j^{\lfloor k_j.(1-\frac{1}{m_j})\rfloor})$$

Now each subset $J_{\ell}$ can be uniquely written as a disjoint union $J_{\ell}=H_{\ell}\cup K_{\ell}$, with $H_{\ell}\subset H:=\{1,2,\dots, (n-p)\}$, and $K_{\ell}\subset K:=\{(n-p+1),\dots, n\}$. Thus $dz_{(J)}:=\otimes _{\ell=1}^{\ell=N}dz_{(J_{\ell})}=\otimes_{r=0}^{r=q}dz_{(J)}(r)$, where $dz_{(J)}(r):=\otimes_{\{\ell\vert r(\ell)=r\}}dz_{(J_{\ell})}$, and, for each $\ell$, $r_{\ell}$ is defined to be the cardinality of $H_{\ell}$. This expression is the local splitting of the filtration of $Sym^N(\Omega^q(X))_{\vert X_y}$ described above.

This filtration naturally induces a corresponding filtration on $S^N(\Omega^q(X\vert \Delta))_{\vert X_y}$, with graded pieces locally given in the coordinates above by: 
$$\frac{dz_{(J)}}{z_{(J)}}:=\frac{1}{z_{(J)}}.\otimes_{r=0}^{r=q}dz_{(J)}(r)=\otimes_{r=0}^{r=q}\frac{dz_{(J)}(r)}{z_{(J)}(r)},$$ where, as above: 
$$z_{(J)}:=\Pi_{j=1}^{j=(n-p)}(z_j^{\lfloor k_j.(1-\frac{1}{m_j})\rfloor}),$$ while:
$$z_{(J)}(r):=\Pi_{j=1}^{j=(n-p)}(z_j^{\lfloor k_j(r).(1-\frac{1}{m_j})\rfloor}),$$ the integers $k_j,$ and $k_j(r)$, if $r>0$, being defined as follows:

1. $k_j$ is the number of $\ell$'s in $\{1,\dots,N\}$  such that $j\in J_{\ell},$

2. $k_j(r)$ is the number of $\ell$'s in $\{1,\dots,N\}$  such that $j\in J_{\ell}$ and $r(\ell)=r,$ for $r=0,\dots,q$.

\

From the definitions and equality above follows that, for any $j=1,\dots, (n-p)$:

3. $z_{(J)}(0):=\Pi_{j=1}^{j=(n-p)}(z_j^{k'_j})$, with: $k'_j:=\lfloor k_j.(1-\frac{1}{m_j})\rfloor-\sum_{r=1}^{r=q}\lfloor k_j(r).(1-\frac{1}{m_j})\rfloor$.

\

The following trivial estimate will be crucial, here:

\begin{lemma}\label{floor} $\lfloor k_j(0).(1-\frac{1}{m_j})\rfloor\leq k'_j\leq q+\lfloor k_j(0).(1-\frac{1}{m_j})\rfloor$, for any $N,n,p$ and $j\in \{1,\dots,(n-p)$.
\end{lemma}

{\bf Proof:} This follows from the fact (applied to $x_r:= k_j(r).(1-\frac{1}{m_j}))$, that the integral part of the sum of $q$ nonnegative real numbers $x_r, r=1,\dots, q,$ lies between the sum of their integral parts, and this same sum increased by $q$ $\square$

\

$\bullet$ The second step of the proof of \ref{f*surj} consists in showing that any nonzero section $s$ of $S^N(\Omega^q(X\vert \Delta))$, when restricted to $X_y$, has its image contained in the smallest piece of the above filtration, described as $\frac{dz_{(J)}(0)}{z_{(J)}(0)}$ in the local coordinates above. But this follows immediately from the hypothesis made, and the filtration given, since otherwise the graded pieces $\frac{dz_{(J)}(r)}{z_{(J)}(r)}$, for $r>0$, would induce non-zero sections over $X_y$ of some sheaf $\otimes_{s=0}^{s=q} S^{N_s}(\Omega^{s}(X\vert \Delta))_{\vert X_y}\otimes \Lambda^{q-s+1}(M)$ , which are supposed not to exist.

$\bullet$ The fird step will show now that the given section is of the form $f^*(w)$ for some local section $w$ of $Sym^N(\Omega^q_Y)$ near $y\in Y$, generic. From the description of the last term $\frac{dz_{(J)}(0)}{z_{(J)}(0)}$ of the above filtration, we see from the last lemma \ref{floor} above, that $s$ is, near $y\in Y$, a section of $f^*(Sym^N(\Omega^q_Y))(q.D)$, where $D$ is the support of $\Delta$. Notice now that the above filtration is compatible with tensor products. Thus replacing $s$ with a sufficiently high tensor power $s^{\otimes k}$, the (obvious) lemma \ref{poles} below shows that $s$ has, in fact, no poles along $D$. Thus $s=f^*(w)$ for some section $w$ of $f^*(Sym^N(\Omega^q_Y))$ near $y$, generic in $Y$.

\begin{lemma}\label{poles} Let $s$ be a meromorphic section of a line bundle $L$ on the complex manifold $X$. Let $D$ be a reduced divisor on $X$. Assume that, for for some given integers $k>q>0$, $s^{\otimes k}$ has poles of order at most $q$ along $D$, and is holomorphic outside of $D$. Then $s$ is holomorphic.
\end{lemma}

$\bullet$ The fourth and last step consists in showing that the (generically uniquely) defined section $w$ of $Sym^N(\Omega^q_Y)$ extends meromorphically to $Y$ as a section of $S^N(\Omega^q(Y\vert \Delta_Y))$, where $\Delta_Y:=\Delta_{(f\vert \Delta)}$. Because the sheaf $Sym^N(\Omega^q_Y)$ is locally free, it is sufficient, by Hartog's theorem, to show this extension at the generic point of any prime divisor $E$ on $Y$. If $f$ is submersive over the generic point $y$ of $E$, and no component of $D$, the support of $\Delta$, is mapped onto $E$, this is clear, since one can take a local section $g$ of $f$ over a neighborhood $U$ of $y$ in $Y$, restrict $s$ to the image of this section. Then $g^*(s)$ and $w$ coincide at the generic point of $U$, and so $g^*(s)$ is the desired extension of $w$. We shall adapt this easy argument to the orbifold context. Let indeed, $m_E$ be the multiplicity of $E$ in the orbifold $\Delta_Y$. There is a irreducible divisor $F\subset X$ such that $m_E=t.m_F$, where $m_F$ is the multiplicity of $F$ of $F$ in $\Delta=\Delta_X$, and $t$ is the positive integer such that $f^*(E)=t.F+\dots$, that is the multiplicity of $F$ in $f^*(E)$. Assume first that $t=1$. Then we have $m_E=m_F$, and moreover a local section $g$ of $f$ exists. The very same argument as above applies, and thus gives the desired extension $g^*(s)$ of $w$, but this time $g^*(s)$ is only a section of $Sym^N(\Omega^q(Y\vert \Delta_Y))$, since $s$ is a section of the restriction to $g(U)$ of $Sym^N(\Omega^q(X\vert \Delta))$, and $g:(Y\vert \Delta_Y)_{\vert U}(X\vert\Delta)_{\vert g(U)}$ is a local isomorphism of orbifolds, by sufficiently shrinking $U$, and choosing $g(y)$ sufficiently generic in $F_y$. The general case where $t> 1$ is now deduced from the case $t=1$ by making a local base change over $h:U'\to U$ which ramifies at order $t$ along $E$. One gets by the preceding argument an extension $\sigma':=(g')^*(s')$ of $w':=h^*(w)$ to $U'$, as a section of $Sym^N(\Omega^q(Y'\vert \Delta_{Y'}))$ over $U'$, where $\Delta_{Y'}$ is defined over $U'$ as the single divisor $E'=h^{-1}(E)$, equipped with the multiplicity $\frac{m_E}{t}$. A simple local computation (in dimension one, in fact) shows that $\sigma'=h^*(\sigma)$, for $\sigma$ a section of $Sym^N(\Omega^q(Y\vert\Delta_{_Y}))$ over $U$. Indeed, the potential poles of $\sigma$ are contained in $E$, with multiplicities at most $t$ (coming from a difference in round-downs), independent on $N$. Thus lemma \ref{poles} above applies again, by replacing $s$ with $s^{\otimes (t+1)}$ $\square$

\section{\bf Orbifold base of a fibration.}\label{ob}

\subsection{Orbifold base}

\

\

In the sequel, $(f\vert \Delta):(X\vert \Delta)\to Y$ will denote a smmoth orbifold $(X\vert \Delta)$, with $X$ projective and $n$-dimensional, together with a fibration $f:X\to Y$, $Y$ being $p$-dimensional.

\begin{definition}\label{dob}
Let $D \subset Y$ be a prime divisor. Then $f^* D = \sum_k t_k. E_k + R,$
where $f(E_k)=D$ and $\codim_Y (f(R)) \geq 2$, i.e. $R$ is an $f$-exceptional divisor.

We define $m(f\vert \Delta;D) := \inf_k \{ t_k.m_{\Delta}(E_k) \}$
to be the multiplicity of the $(f\vert \Delta)$-fibre over a general point of $D$.

Set $\Delta_{(f\vert\Delta)} := \sum_{D \subset Y} (1-\frac{1}{m(f;D)}) D$. The `orbifold base' of $(f\vert \Delta)$ is defined to be $(Y\vert\Delta_{(f\vert \Delta)})$.
When $\Delta=0$, we simply write: $m(f;D)$ and $\Delta_f$.

\end{definition}

{\bf Remark.} 0. The (simple) reason  for this definition is given in remark \ref{rcomp} below.

1. Since the general $f$-fibre is smooth, so irreducible and reduced, 
it is clear that there are at most finitely many prime divisors $D$ such that $m(f\vert\Delta;D)>1$.
Thus $\Delta_f$ (resp. $\Delta_{(f\vert \Delta)})$ is always an integral orbifold divisor (resp. when so is $\Delta)$.

2. $\Delta(f\vert \Delta')\geq \Delta(f\vert \Delta)$ if $\Delta'\geq \Delta$.

3. In general, $f:(X\vert \Delta)\to (Y\vert \Delta_{(f\vert \Delta)})$ is an orbifold morphism only in codimension one, because the multiplicities on the $f$-exceptional divisors on $X$ are not taken into account in the definition, and may be too small. This is remedied (to a certain extent, sufficient for the purposes) by the consideration of `neat models' in the next section.

4. The `classical' (or `divisible') notion {\it does not coincide} with our multiplicity above. With this notation,
it is defined by $ m^*(f\vert \Delta;D):= gcd_k \{ t_k. m_{\Delta}(E_k)\}.$
Thus $m^*(f\vert \Delta,D)$ divides $m(f\vert \Delta,D)$, and $(f\vert \Delta)$
 will have \lq more\rq\ 
multiple fibres than in the classical sense, in general. 

5. There are several main reasons for considering the $inf$, instead of the $gcd$ multiplicities: compatibility with differentials (see \ref{mor}), canonical dimension of orbifold bases (see \ref{L_f}). Also for hyperbolicity reasons, the one can expect the Kobayashi pseudometric of $(X\vert \Delta)$ to be the lift from the one on the orbifold base of its `core' (see \S 7) only with the $inf$ multiplicities.

\begin{remark}\label{rcomp} There is a composition rule\footnote{On suitable `neat' birational models of $f,g$ at least. See below for this notion.} for the base orbifold of the composition $g\circ f$ two fibrations $X \stackrel{f}{\rightarrow} Y \stackrel{g}{\rightarrow} Z$, namely: $\Delta_{g \circ f} = \Delta_{(g, \Delta_{f})}$, when $\Delta=0$. This equality suggests the definition of the base orbifold of $f:(X\vert \Delta_X)\to Y$, given above. We refer to \cite[ch.1.6]{Ca04} and \cite[chap. 3]{Ca07} for more details.\end{remark}

\subsection{Birational (non) invariance}

\begin{definition}
Let \holom{f}{X}{Y} and \holom{f'}{X'}{Y'} be two fibrations, and $\Delta, \Delta'$ be orbifold divisors on $X,X'$ respectively. Assume that $u: (X'\vert \Delta')\to (X\vert \Delta)$ is an `elementary birational orbifold morphism' between smooth orbifolds. 

Then $(f\vert \Delta)$ is said to be `elementarily birationally
equivalent' to $(f'\vert \Delta')$, if there exists a birational map  \holom{v}{Y'}{Y}
making commutative the diagram:
\[
\xymatrix{ 
X'\ar[d]^{f'} \ar[r]^{u} & X
\ar[d]^{f}
\\
Y' \ar[r]^{v} & Y'. }
\]
We say that  $(f\vert \Delta)$ is `birationally
equivalent' to $(f'\vert \Delta')$, and  write $(f\vert \Delta) \sim (f'\vert \Delta')$ if they are connected by an chain of elementarily equivalent fibrations $(f_j\vert \Delta_j)$.
\end{definition}

\begin{lemma}
If $(f\vert \Delta)$ is said to be `elementarily birationally
equivalent' to $(f'\vert \Delta')$, we have the following properties:
\begin{itemize}
\item $v_* (\Delta(f'\vert \Delta'))=\Delta(f\vert \Delta)$,
\item $\kappa(Y'\vert\Delta(f'\vert \Delta')) \leq \kappa(Y\vert\Delta(f\vert \Delta))$,
\item and $\kappa(Y'\vert\Delta(f'\vert \Delta')) = \kappa(Y\vert\Delta(f\vert \Delta))$, if $\kappa(Y\vert \Delta(f\vert \Delta)) \geq 0$.
\end{itemize}
\end{lemma}

{\bf Remark.} The inequality $\kappa(Y'\vert\Delta(f')) \leq \kappa(Y\vert\Delta(f))$ can be strict
if $\kappa(Y\vert \Delta(f\vert \Delta))=-\infty$.  
In particular the canonical dimension of the orbifold base is
{\it not} a birational invariant of fibrations. The next definition remedies this situation, but using, in a first step, a less computable invariant. This drawback will be cancelled using `neat fibrations' below.

\begin{definition}
Let \holom{f}{X}{Y} be a fibration, then we define the canonical dimension of $(f\vert \Delta)$ as $\kappa(f\vert \Delta) := \inf_{(f'\vert \Delta') \sim (f\vert \Delta)} \{\kappa(Y', K_{Y'}+\Delta(f'\vert \Delta')) \}.$
\end{definition}

The canonical dimension $\kappa(f\vert \Delta)$ is now a birational invariant and we can extend the definition
to any rational map \merom{f}{X}{Y} $Y$ arbitrarily singular by resolving the singularities
of $Y$ and the indeterminacies of $f$. We shall now compute this  canonical dimension differently in the next section \ref{diff}.

\begin{remark}\label{rbir} The study of orbifold birational maps needs to be completed on several essential points. In particular:

1. They should be extended to the category of log-canonical and klt orbifolds.

2. It is not known whether the orbifold bases of `neat' fibrations are birational invariants (in the orbifold category), and this even for the three central fibrations of the decompositions described below (the Moishezon-Iitaka fibration, the `$\kappa$-rational quotient', and the `core', see remarks \ref{rM}, \ref{rr} and \ref{rbc})
\end{remark}

\subsection{The differential sheaf of a fibration}\label{diff}

\begin{definition}
Let $E$ be a locally free sheaf on a complex manifold $X$. Let $\sF \subset E$
be a coherent subsheaf of $E$. The saturation $\sF^{\sat}$ is the kernel
of the map:  $E \rightarrow  (E/\sF)/(\Tor E/\sF)$. It is also the largest subsheaf of $E$ containing $\sF$, and identical with $\sF$ at the generic point of $X$.
\end{definition}

\begin{definition} Let $(X\vert \Delta)$ be a smooth orbifold.
Let \holom{f}{X}{Y} be a holomorphic fibration between projective manifolds, with $p:=\dim Y$. 
 For all $N \geq 1$, the 
canonical injection $0 \rightarrow f^* \Omega^p_Y \rightarrow \Omega^p_X$
induces a injective morphism $f^* (N K_Y)=f^* ((\Omega_Y^p)^{\otimes N}) \rightarrow \Sym^N (\Omega_X^p)$
and we define: $L^N_{(f\vert \Delta)} := (f^* (N K_Y))^{\sat}$ 
to be the saturation of $f^*(N K_Y)$ in $S^N\Omega^p(X\vert \Delta)$.
\end{definition}

By definition $L^N_{(f\vert \Delta)}$ is a coherent rank $1$ subsheaf of $S^N\Omega^p(X\vert \Delta))$ and
 it is easily seen to be a birational invariant of the fibration $(f\vert \Delta)$: elementary birational equivalences of fibrations indeed induce isomorphisms at generic points, but not on global sections in general. Conditions for this are given in theorem \ref{L_f} below. Moreover, for $N$ sufficiently large and divisible, 
 $f^* (N (K_Y+\Delta(f\vert \Delta)))\subset L^N_{(f\vert \Delta)}$ outside the $f$-exceptional divisors of $X$. 
 
 This is an elementary local computation that we explain, to simplify notations, when $\Delta=0$, $n=dim X=2$ and $p=dimY=1$. The map $f$ is then given in suitable local coordinates $(u,v)$ on $X$, and $y$ on $Y$, by: $f(u,v)=y=u^k$, near a generic smooth point of a component of multiplicity $k$ of some fibre $X_0=f^{-1}(0)$ of $f$. Let $m$ be the $inf$ multiplicity of this fiber, so that $k\geq m$. Thus $f^*(Nm.(K_Y+\Delta(f))$ is generated, as an $\sO_X$-module, by: $f^*((\frac{dy}{y^{(1-\frac{1}{m})}})^{\otimes N.m})=k^{N.m}.u^{(k-m).N}.(du)^{\otimes N.m}$. This gives the said inclusion, and equality when $k=m$, in this case.

\begin{definition}\label{neat}
A fibration \holom{f}{X}{Y} is `neat' if there exists a birational morphism
\holom{u_0}{X}{X_0} such that $X_0$ is smooth and if every $f$-exceptional prime divisor $E \subset X$  is also $u_0$-exceptional. 
\end{definition}

Neat fibrations are constructed as follows: for any fibration
 \holom{f}{X}{Y} with $X$ smooth, let $Y' \rightarrow Y$ be a birational morphism such that 
$X \times_Y Y' \rightarrow Y'$ is flat (which exists, by Raynaud's flattening theorem).
A desingularisation $X' \rightarrow X \times_Y Y'$ of the
main component of $X \rightarrow X \times_Y Y'$, induces a neat fibration \holom{f'}{X'}{Y'}, 
by taking for $u_0$ the natural projection $u_0:X'\to X$.

\begin{theorem}\label{L_f}
Let  $(X\vert \Delta)$ be a smooth orbifold, and \holom{f}{X}{Y} be a neat fibration, then:
$H^0(X, L_{(f\vert \Delta)}^N) \simeq f^* (H^0(Y, N (K_Y + \Delta_{(f\vert \Delta)})))$, for all $N$ large and sufficiently divisible. 

For any rational fibration $f:X\dasharrow Y$, we can now define: $\kappa(f\vert \Delta):=
\kappa(Y'\vert\Delta_{(f'\vert \Delta')})$, if $f':X'\to Y'$ is a neat fibration birationally equivalent to $f$, with $u:(X'\vert \Delta')\to (X\vert \Delta)$ an elementary birational orbifold morphism. The R.H.S of the equality being independent on the choice of $(X'\vert \Delta')$ by theorem 3.5 and the remark following definition 3.6.\end{theorem}

{\bf Idea of proof:} First assertion: The fibration being neat, Hartog's theorem shows that we do not change the space of sections on the left-hand side by allowing poles of arbitrary orders on sections of $L^N_{(f\vert \Delta)}$ along divisors of $X$ which are $f$-exceptional, hence $u_0$-exceptional. We thus assume there are no such divisors. 

By the definition of $L^N_{(f\vert \Delta)}$, we have then a natural inclusion: $N (K_Y + \Delta(f))\subset L_{(f\vert \Delta)}^N$ with cokernel supported on a divisor $D$ of $X$ `partially supported on fibres of $f$', which means that if $E\subset Y$ is an irreducible component of $f(D)$, then $f^{-1}(E)$ has an irreducible component $D'$ not contained in $D$, and such that $f(D')=E$.

An elementary general lemma now shows that, in this situation, the global sections of the two said sheaves coincide.

Second assertion: when $f$ is neat, it is an immediate consequence of the first assertion. The general case follows from the fact that both sides are birational invariants of the fibration $\square$

\begin{remark} By contrast, we do not know if the base orbifolds of two birationally equivalent neat fibrations are birationally equivalent. This is an important question, seemingly quite difficult.
\end{remark}

\subsection{Fibrations of base-general type}

\

\begin{definition}
A meromorphic fibration \merom{(f\vert \Delta)}{(X\vert \Delta)}{Y} is of base-general type if $\kappa(f\vert \Delta)=\dim Y>0$, with $(X\vert \Delta)$ a smooth orbifold. (We shall abbreviate `base-general type' to `general type').
\end{definition}

{\bf Notation.} We denote by $\sF(X\vert \Delta)$ the set of fibrations of (base-)general type on $(X\vert \Delta)$ 
(up to birational equivalence).

\begin{definition}\label{dkL}
Let $X(X\vert \Delta)$ be a smooth, projective orbifold. A `Bogomolov sheaf' on $(X\vert \Delta)$ is a rank one 
coherent, saturated, subsheaf $L \subset \Omega_X^p$ such that: $\kappa(X\vert \Delta,L) = p>0$.

By definition, $\kappa(X\vert \Delta,L):=\overline{lim}_{N>0}\frac{log(h^0(X,(L^{N,\Delta}))}{log(N)}$, 
where $L^{N,\Delta}$ is the saturation in $S^N\Omega^p(X\vert \Delta)$ of the image of $f^*(N.K_Y)$.
\end{definition}

{\bf Notation.} We denote by $Bog(X\vert \Delta)$ the set of Bogomolov sheaves on $X$. 
Bogomolov sheaves are thus the rank one subsheaves of $\Omega_X^p$ of maximal $\Delta$-positivity, by
 the next theorem, essentially due to Bogomolov. It rests on Deligne's closedness of logarithmic forms on $(X\vert \lceil \Delta\rceil)$.

\begin{theorem}\label{bog}(\cite{Bo78}, \cite{Ca07})
Let $(X\vert \Delta)$ be a smooth projective orbifold and $L \subset \Omega_X^p$ a rank 1 subsheaf. Then: $\kappa(X\vert \Delta,L) \leq p$, 
and if equality holds, there exists a meromorphic fibration \merom{f}{X}{Y} such that
$L= f^* \Omega_Y^p$ holds above the generic point of $Y$.
\end{theorem}

Together with theorem \ref{L_f}, Bogomolov's theorem leads immediately to the following geometric description of Bogomolov sheaves:

\begin{theorem}
\label{theoremmapsbogomolovsheaves}\cite[Th\'eor\`eme 8.9]{Ca07}
The map $\holom{L}{\sF(X\vert \Delta)}{Bog(X\vert \Delta)}, [f] \mapsto L_f$ is bijective.
\end{theorem}

{\bf Remark.} It is essential for the theorem to use `inf-multiplicities', instead of `gcd-multiplicities'. The theorem \ref{theoremmapsbogomolovsheaves} interprets geometrically the ``Bogomolov sheaves", which were only partially interpreted in \ref{bog}.

An important property of fibrations of general type is:

\begin{theorem}\label{fgtah} \cite[Th\'eor\`emes 8.17-19]{Ca07} Let $(f\vert \Delta):(X\vert \Delta)\dasharrow Y$ be a fibration, with $(X\vert \Delta)$ smooth projective. Then $f:X\dasharrow Y$ is `almost holomorphic' if either $\kappa(f\vert \Delta)=dim Y$ is of general type, or if $\kappa(f\vert \Delta)\geq 0$ and $(X\vert \Delta)$ is finite.

Recall that $f$ is `almost holomorphic' if its indeterminacy locus does not meet its general Chow-theoretic fibre.
\end{theorem}

\begin{example}\label{ernh} When $(X\vert \Delta)$ is logarithmic, a rational map $(f\vert \Delta):(X\vert \Delta)\dasharrow Y$ with $\kappa(f\vert \Delta)\geq 0$ need not be almost holomorphic. Take $X=\Bbb P^2$, $\Delta$ the union of two lines meeting in a point $a$, and $f:X\to \Bbb P^1$ the linear projection from $a$. 
\end{example}

\subsection{Orbifold fibres and suborbifolds.} \label{of}

\

\

Let $(f\vert \Delta):(X\vert \Delta)\to Y$ be a fibration from the smooth projective orbifold $(X\vert \Delta)$. For $y\in Y$, let $\Delta_y$ be the restriction of the divisor $\Delta$ to the fibre $X_y:=f^{-1}(y)$ of $f$. From Sard theorem, it follows that , for $y\in Y$ generic, $\Delta_y$ is an orbifold divisor of $X_y$, and that the orbifold $(X_y\vert \Delta_y)$ is smooth. We call it the generic orbifold fibre of $(f\vert \Delta)$.

Moreover, if $(f'\vert \Delta'):(X'\vert \Delta')\to Y'$ is a fibration birationally equivalent to $(f\vert \Delta):(X\vert \Delta)\to Y$, it is easy to check that for some dense Zariski open subset $U=U'$ in $Y\sim Y'$, the orbifold fibres of $(f\vert \Delta)$ and $(f'\vert \Delta')$ over $y$ are birationally equivalent (in the orbifold sense).

\begin{remark}\label{eso} When $(f\vert \Delta):(X\vert \Delta)\dasharrow Y$ is only rational, the generic orbifold fibre of $(f\vert \Delta)$ is not well-defined, up to birational orbifold equivalence (unless $f$ is almost holomorphic). Even its canonical dimension is not well-defined up to birational equivalence. Consider $\Bbb P^2$, and the map $g:X\dasharrow \Bbb P^1$ defined by a pencil of conics through $4$ points. Let $u:X\to \Bbb P^2$ be the blow-up in these $4$ points, with exceptional divisor $E$ consisting thus of four $(-1)$-curves. Let $\Delta$ be the orbifold divisor on $X$ obtained in attributing the multiplicity $m$ to each component of $E$. The family of the strict transforms of the smooth members of the initial pencil of conics are fibres of $f:=g\circ u$, and their orbifold canonical dimension is thus $-\infty$ (resp. $0$, $1)$ if $m=1$ (resp. $m=2$, $m\geq 3)$.
\end{remark}

\begin{definition} Let $(X\vert \Delta)$ be a smooth projective orbifold, and $j:V\to X$ the inclusion of an (irreducible) subvariety not contained in $Supp(\Delta)$. We define a `restriction' of $\Delta$ to $V$ as any smooth projective orbifold pair $(W\vert \Delta_W)$ together with a birational map $g:W\to V$ such that the composed map $j\circ g:(W\vert \Delta_W)\to (X\vert \Delta)$ is an orbifold morphism.  

For any given birational map $g:W\to V$, with exceptional divisor $E$, such that $(E\cup(j\circ g)^{-1}(Supp(\Delta)))$ is a divisor of normal crossings, there exists a smallest restriction $(W\vert \Delta_W)_{min}$ of $\Delta$ to $V$ defined on $W$. \end{definition}

\begin{remark}

1. However, if $h:W'\to W$ is another birational morphism such that $g\circ h:W'\to V$ satisfies the normal crossings condition above, it is not known presently if the minimal restriction $(W'\vert \Delta')_{min}$ of $\Delta$ to $W'$ is an orbifold morphism\footnote{Except when $V$ is a curve, a case considered in more details in \S \ref{orc}}. If this were so, the notion of minimal restriction would be birational (for given $(X\vert \Delta)$, $V)$.

2. However, the example in remark \ref{eso} shows that, even for curves, and for generic fibres of a  rational fibration, there is no possibility to define a birationally well-defined notion of minimal restriction which is independant on the birational model of $(X\vert \Delta)$ on the corresponding strict transform of $V$. (Except for the generic member of a `base-point free' covering family of suvarieties).

For this reason, we give the following definition:
\end{remark}

\begin{definition}\label{dso} Let Let $(X\vert \Delta)$ be a smooth projective orbifold, and $V\subset X$ a subvariety not contained in $Supp(\Delta)$. We say that the restriction of $\Delta$ to $V$ has a (birationally invariant) property $P$ if there is some (smooth) orbifold-birational model $(X'\vert \Delta')$ of $(X\vert \Delta)$ such that the (minimal) restriction of $\Delta'$ to the strict transform $V'$ of $V$ on $X'$ has the property $P$. 

The properties $P$ which will be important here will be: being `special', or with $\kappa=0$, or with $\kappa_+=-\infty$. See below.
\end{definition}

\begin{proposition} Let $V_t$ be the generic member of a covering family of suvarieties of $X$. Having property $P$ is independent on the birational model of the (smooth and projective) $(X\vert \Delta)$. The `genericity' here depends on the model. (This applies, in particular, to orbifold fibres of rational fibrations).
\end{proposition}

\

\section{\bf Orbifold additivity}

The orbifold version of Iitaka's $C_{n,m}$-conjecture is:

\begin{conjecture}($C_{n,m}^{\orb}$)\label{cnmorb}  Let $(f\vert \Delta):(X\vert \Delta)\to Z$ be a fibration, with $(X\vert \Delta)$ smooth and $X$ projective. Then: $\kappa(X\vert\Delta) \geq \kappa(Y_z\vert\Delta_z)+ \kappa(g\vert \Delta).$ Here $\Delta_z$ is simply the restriction of $\Delta$ to $X_z:=g^{-1}(z),z\in Z$ generic.
\end{conjecture}

Our main technical result is the following:

\begin{theorem}\label{vieh}
The $C_{n,m}^{\orb}$-conjecture is true if $\kappa(f\vert\Delta) = \dim Z$ (ie: if $(f\vert \Delta)$ is of general type).

In this case: $\kappa(X\vert \Delta)=\kappa(X_z\vert \Delta_z)+dim Z$.
\end{theorem}

\begin{corollary} \label{addcomp} Let $(X\vert \Delta)$ be a smooth projective orbifold, and let \holom{f}{X}{Y} and \holom{g}{Y}{Z} be fibrations.  
Then:  $\kappa(f\vert \Delta) = \kappa(f_z\vert \Delta_z) + \dim Z$,  if $g \circ f$ is of general type, denoting $(f_z\vert \Delta_z):(X_z\vert \Delta_z)\to Y_z$ the restriction over a general $z\in Z$.

Thus: if $(g \circ f\vert \Delta)$ and $(f_z\vert \Delta_z)$ are of general type, then $(f\vert \Delta)$ is of general type.\end{corollary}

This is a fairly direct application of the following semi-positivity result ([Ca04, theorem 4.11, p. 567], differently formulated), extending \cite{V}:

\begin{theorem}\label{vieh'}\cite[Theorem 4.11, p. 567]{Ca04} Let $(X\vert \Delta)$ be a smooth geometric, with $X$ projective. Let $f:X\to Z$ be a fibration, with $Z$ smooth. Assume that $f$ is `prepared' (see \cite[1.1.3, p. 508]{Ca04})\footnote{This can be achieved by an elementary orbifold modification of $(X\vert \Delta)$, which does not change the space of sections of $N.(K_X+\Delta)$. The term is introduced in \cite{V}.}.

Let $(Z\vert \Delta_Z)$, with  $\Delta_Z:=\Delta_{(f\vert \Delta)}$ be the orbifold base of $(f\vert \Delta)$.  Define 
$K_{X/(Z\vert \Delta_Z)}:=K_X-f^*(K_Z+\Delta_Z)$.

Then , $f_*(N.(K_{X/(Z\vert \Delta_Z)}+\Delta))$ is weakly-positive in Viehweg's sense for any $N>0$ such that $N.(K_{X/(Z\vert \Delta_Z)}+\Delta)$ is Cartier.

In particular, the restriction of $\sF:=m.(N.(K_{X/(Z\vert \Delta_Z)}+\Delta))+f^*(B)$ over the generic fiber $X_y$ of $f$  is generated by the global sections of $\sF$ on $X$, for any given $N$ as above, and any big $\Q$-divisor $B$ on $Y$, if $m$ is sufficiently large (depending on $N,B)$.

\end{theorem}


A result of a similar nature, more precise (the conclusion is nefness), is shown in \cite{ka}, when the orbifold fibres have a trivial canonical bundle,.

\

The proofs are essentially adaptations to the orbifold context of Viehweg's proof of the $C_{n,m}$-conjecture
for fibrations with base a variety of general type (\cite{V}). We shall not give the rather technical and lengthy details here. 
The range of applications
is however considerably increased by the orbifold context. It is, for example, one of the crucial ingredients used in \cite{MH} and \cite{MH'}. A different proof is given in \cite{BP}.

\begin{corollary}
Let \holom{(f\vert \Delta)}{(X\vert \Delta)}{Z} be a fibration of general type, then: $\kappa(X\vert \Delta) \geq \kappa(X_z\vert \Delta_{z}) + \dim Z.$

\end{corollary}

\begin{corollary}
\label{ck0}
Let $(X\vert \Delta)$ be a smooth projective orbifold with $\kappa(X\vert \Delta)=0$. There is no fibration $f:Y\dasharrow Z$ of general type on $(X\vert \Delta)$.
\end{corollary}

\begin{example}
\label{ctorics}
Let $(X\vert \Delta)$ be the smooth projective orbifold with $X$ a projective {\it toric} manifold, and $\Delta$ the reduced anticanonical divisor complement of the open orbit. There is no fibration $f:X\dasharrow Z$ of general type on $(X\vert \Delta)$. (Example suggested by a question of M. Musta\c ta).
\end{example}

\begin{corollary}
\label{cfanos}
Let $(X\vert \Delta)$ be a smooth projective Fano orbifold (ie: such that $-(K_X+\Delta)=H$ is ample on $X)$. There is no fibration $f:X\dasharrow Z$ of general type on $(X\vert \Delta)$.
\end{corollary}

{\bf Proof:} Choose a reduced $D\in \vert N.H\vert$, for $N$ large and divisible such that $\Delta':=\frac{1}{N}.D+ \Delta$ has a support of normal crossings. Then $(X\vert \Delta')$ is smooth, and has trivial canonical bundle. Since there is, by \ref{ck0}, no fibration of general type on $(X\vert \Delta')$, there is, a fortiori, no fibration of general type on $(X\vert \Delta)$ $\square$

\begin{remark} The conjecture $C_{n,m}^{orb}$ implies much more: if $(f\vert \Delta):(X\vert \Delta)\dasharrow Z$ is a fibration from $(X\vert \Delta)$, smooth and projective, it should hold that: 

1. $\kappa(f\vert \Delta)\leq 0$ if $\kappa(X\vert \Delta)= 0.$

2. $\kappa(f\vert \Delta)=-\infty$ if $\kappa(X\vert \Delta)=-\infty$ and $\kappa(X_z\vert \Delta_z)\geq 0.$

In fact, even `more' should be true.
 \end{remark}

This `more' is best formulated with the obvious orbifold versions of $\kappa_{++}$ and $\kappa_+$:

\begin{definition}\label{dok++} Let $(X\vert \Delta)$ be a smooth projective complex connected orbifold. Let: $\kappa_{++}(X\vert \Delta):=max_{\{L,p>0\}}\{\kappa(X\vert \Delta,L)\}$, where $L\subset \Omega_X^p$ is a  rank $1$ coherent subsheaf. (See definition \ref{dkL} for the definition of $\kappa(X\vert \Delta,L)).$

Define also: $\kappa_{+}(X\vert \Delta):=max_{\{f\}}\{\kappa(f\vert \Delta)\}$, where $(f\vert \Delta):(X\vert \Delta)\dasharrow Z$ ranges over all fibrations  defined on $X$.

Thus: $dim(X)\geq \kappa_{++}(X\vert \Delta)\geq \kappa_{+}(X\vert \Delta)\geq \kappa_{}(X\vert \Delta)$ for any $(X\vert \Delta)$.
\end{definition}

\begin{remark}

1. Notice that, when $\Delta=0$, it is not obvious that the  old $\kappa_+(X)$ of \ref{dk+} coincides with the new $\kappa_+(X\vert 0)$ in the sense of \ref{dok++}, which is greater. But conjectures \ref{cnmorb} above and \ref{ck+} imply that $\kappa_+(X\vert 0)=\kappa_+(X)$.

2. The conjecture \ref{cok+} in \S \ref{orc} asserts that, just as when $\Delta=0$, the condition $\kappa+=-\infty$ is equivalent to orbifold rational connectedness. \end{remark}

\begin{conjecture} \label{cok++} If $\kappa(X\vert \Delta)\geq 0$, then $\kappa_{++}(X\vert \Delta)=\kappa(X\vert \Delta)$.

If $\kappa_+(X\vert \Delta)=-\infty$, then $\kappa_{++}(X\vert \Delta)=-\infty$.
\end{conjecture}

The conjectures $C_{n,m}^{orb}$ and \ref{cok++} will be deduced from standard conjectures of the LMMP in \S\ref{abund}. The central conjecture (see \ref{pfcm}) concerning families of canonically polarized manifolds will be shown to follow from conjectures $C_{n,m}^{orb}$ and \ref{cok++} in \S\ref{fcpm}, and thus from standard conjectures of the LMMP.

\section{\bf Special Orbifolds}

\subsection{Definition and Main examples.}

\begin{definition}
A smooth projective orbifold $(X\vert \Delta)$ is `special'
if there does not exist a dominant rational map \merom{(f\vert \Delta)}{(X\vert \Delta)}{Y} 
of base-general type (ie: with $\kappa(f\vert \Delta)=\dim Y>0)$.

A variety $X$ is special if so is $(X'\vert 0)$, for some (or any) smooth model $X'$ of $X$.

A quasi-projective manifold $U$ is special if so is $(X\vert D)$ for some (or any) smooth compactification $X$ of $U$ by a normal crossing divisor $D$.
\end{definition}

{\bf Remarks.}

0. Being special is indeed a birational property.

1. A (smooth projective, as always in the sequel) orbifold $(X\vert \Delta)$ is special if and only
if there does not exist a Bogomolov sheaf on $(X\vert \Delta)$ (by theorem  \ref{theoremmapsbogomolovsheaves}).

2. An orbifold of general type and positive dimension, is not special (consider the identity map).

3. An orbifold  curve $(X\vert \Delta)$ is special if and only $deg(K_X+\Delta)\leq 0$, since a fibration $f:X\dasharrow Y$  is either
the identity map or the constant map.

4. A rationally connected manifold $X$ is special, since $\kappa_{++}(X)=-\infty$.

5. A Fano orbifold is special, by corollary \ref{cfanos}.

6. Any orbifold with $\kappa=0$ is special, by corollary \ref{ck0}.

7. Any orbifold with $\kappa_+=-\infty$ is obviously special.

8. The previous two examples are in fact the `building blocks' of special orbifolds. First, the theorem \ref{tcompspec} shows that towers of fibrations with general orbifold fibres having either $\kappa=0$, or $\kappa_+=-\infty$ are special. Conversely, we shall see in section \ref{decc} that, conditionally in $C_{n,m}^{orb}$, any special orbifold decomposes canonically as a tower of such fibrations.

\begin{theorem}\label{tcompspec} Let $(f\vert \Delta):(X\vert \Delta)\to Y$ a neat fibration, with $(X\vert \Delta)$ smooth projective. Assume that its general orbifold fibre $(X_y\vert \Delta_y)$ and orbifold base $(Y\vert \Delta(f\vert \Delta))$ are special. Then $(X\vert \Delta)$ is special.
\end{theorem}

The proof is sketched in \ref{extspec}.

\begin{remark}

1. This statement is (very) false if one does not take into account the orbifold structures on both fibres and base (see example \ref{exampleiitaka} when $\Delta=0)$.

2. Orbifolds with either $\kappa=0$ or $\kappa_+=-\infty$ are special, while manifolds of
general type are not special, so one might
ask if further relations between $\kappa$ and specialness hold. This is not the case:
for every $n \in \N$ and every $k \in \{ - \infty, 1, \ldots, n-1 \},$
there exist special manifolds $X$ of dimension $n$ and $\kappa=k$. For example, any 
hyperplane section $X \subset \PP^{n+1-k} \times \PP^k$ of bidegree $(n+2-k, d)$, with $d>k+1,$ and containing a section over $\PP^k$, has $dim X=n$, $\kappa(X)=k$ if $k>0$, and is special. 
\end{remark}

\subsection{Criteria for Special Orbifolds.}

\begin{theorem}\label{sconn}
A smooth projective orbifold is special if (and only if) two generic points are joined by a chain of special suborbifolds (ie: images of orbifold morphisms from special smooth orbifolds).
\end{theorem}

{\bf Proof:} It is sketched in the next section \ref{core}. See remark \ref{propc}
 $\square$

\begin{remark}
When $\Delta=0$, \ref{sconn} does not generalize to arbitrary singular varieties (eg. if $X$ is a cone over a variety of general type $Y)$. It should however hold for varieties with log-canonical singularities.

A special manifold does not necessarily admit non-trivial chains of special subvarieties (eg. if $X$ is a simple
abelian variety).
\end{remark}

\begin{theorem}\label{etal}
Let \merom{g}{X'}{X} be a dominant map. 

$\bullet$ If $X'$ is special, $X$ is special.

$\bullet$ If $X$ is special, and if $g$ is regular and \'etale, then $X'$ is special.  
\end{theorem}

The proof of the seemingly easy second assertion requires the difficult result \ref{vieh}.

\begin{theorem}\label{koorb} Let $\varphi:\C^n\to X$ a non-degenerate meromorphic (possibly transcendental) map. Then $X$ is special.
\end{theorem} 

This extends a former result of Kobayashi-Ochiai asserting that $\kappa(X)< \dim X$ under the above hypothesis of non-degeneracy (which means that $\varphi$ is holomorphic and submersive at some point). More general versions of \ref{koorb} are given in \cite{Ca04}.

\subsection{Special Surfaces.} 

\

\

They can be easily described, using classication. Such a simple description however fails in higher dimensions.

\begin{proposition}
The special surfaces $X$ are exactly the following ones (up to birational equivalence):
\begin{itemize}
\item $\kappa(X)=-\infty$: then $X \simeq \PP^1 \times C$ with $g(C)=0$ or $1$.
\item $\kappa(X)=0$: then $X$ is a K3 surface or abelian.
\item $\kappa(X)=1$: (after a suitable \'etale cover of $X$) the Moishezon-Iitaka elliptic fibration $X \rightarrow C$ has either $g(C)=1$
and no multiple fibre, or $g(C)=0$, and at most 2 multiple fibres of coprime multiplicities.
\end{itemize}
\end{proposition}

\begin{definition}\label{ws} A manifold $X$ is `weakly special' if no finite \'etale cover $X'\to X$ has a meromorphic map $f:X'\dasharrow Y$ onto a variety $Y$ of general type with $\dim Y>0$.

 It follows from theorem \ref{etal} that if $X$ is special, it is weakly special, the converse being true when $\dim X\leq 2$. Indeed:
\end{definition}

\begin{corollary}
Let $X$ be a smooth surface. The following are then equivalent:
\begin{itemize}
\item $X$ is special.
\item $X$ is `weakly special'.
\item $\kappa(X)\leq 1$, and $\pi_1(X)$ is almost abelian.
\item There exists a nondegenerate map \merom{h}{\C^2}{X}.
\end{itemize}
\end{corollary}

\begin{remark}\label{bt}

1. None of these characterizations of specialness extends to higher dimensions. 

2. The property of being special
is thus invariant by deformation for surfaces.

3. There exists weakly special threefolds $X's$ which are not special. See example \ref{btex} for a sketch of their construction). These threefolds show that orbifolds are needed in the birational classification theory of projective varieties (multiple fibres cannot always be eliminated by finite \'etale covers), and also permit to test some conjectures in arithmetics and hyperbolicity (see section \ref{conj} below).
\end{remark}

\section{\bf The Core}\label{core}

The main structure result of this text is:

\begin{theorem}\label{c} For any smooth complex projective orbifold $(X\vert \Delta)$, there exists a unique map $c=c_{(X\vert \Delta)}:X\dasharrow C(X\vert \Delta)$, called `the core' of $(X\vert \Delta)$ such that:
\begin{itemize} 
\item The `general' orbifold fibre $(X_c\vert \Delta_c)$ of $c$ is special.
\item $\kappa(c\vert \Delta)=\dim C(X\vert \Delta)\geq 0$.
\end{itemize}

Moreover, $c$ is almost holomorphic (by theorem \ref{fgtah}).
\end{theorem}

\begin{remark}\label{rbc} We do not know if the orbifold bases of `neat' models of $c$ are birational invariants. 
\end{remark}

{\bf Idea of proof:} Uniqueness follows from lemma \ref{u} below. 

Existence: Let $f:X\dasharrow Y$ a fibration with $d:=\dim Y$ maximum such that $\kappa((f\vert \Delta))=\dim Y$. 
Thus $d=0$ (resp. $d=\dim X$) if and only if $(X\vert \Delta)$ is special (resp. of general type). In the general case, we have to show that the fibres of $(f\vert \Delta)$ are special. Assume not. Then by Chow space theory and countable upper semi-continuity of the above dimension $d$, one can construct a factorisation $f=g\circ h$, where: $h:X\dasharrow Z$, $g:Z\to Y$ with $\dim Z>\dim Y$ and such that for general $y\in Y$, the map: $(h_y\vert \Delta_y):(X_y\vert \Delta_y)\to Z_y$ is of general type. But it then follows from corollary \ref{addcomp} that $(h\vert \Delta):(X\vert \Delta)\to Z$ is a fibration of general type, contradicting the maximality of $\dim Y$ $\square$

\begin{lemma}\label{u} Let $(h\vert \Delta):(X\vert \Delta)\dasharrow Z$ be a fibration with general orbifold fibre special (resp. let $(g\vert \Delta):(X\vert \Delta)\dasharrow Y$ be of general type). Then there exists a unique map $u:Z\to Y$ such that: $g=u\circ h$.
\end{lemma}

{\bf Idea of proof:} Otherwise the images of the general fibres $X_z$ of $h$ by $g$ are positive-dimensional. But one can show that the restriction of the corresponding Bogomolov sheaf $L_{(g\vert \Delta)}$ on $(X\vert \Delta)$ to such an orbifold fibre $(X_z\vert \Delta_z)$ were then a non-zero Bogomolov sheaf on $(X_z\vert \Delta_z)$, which is special. A contradiction $\square$

\begin{remark}\label{propc} The map $c$ enjoys the following properties (as may be seen from the above sketch of proof of \ref{c}):

1. $\dim C(X\vert \Delta)=n$ (resp. $0)$ iff  $(X\vert \Delta)$ is of general type (resp. special).

2. If $V\subset X$ is a subvariety which meets a (suitable) general fibre $X_c$ of $c_X$, and is the image of an orbifold morphism $(W\vert \Delta_W)\to (X\vert \Delta)$, with $(W\vert \Delta_W)$ special, then $V\subset X_c$. This property implies theorem \ref{sconn}.

3. If $u:(X'\vert \Delta')\to X$ is a surjective orbifold morphism, there exists a unique map $c_u:C(X'\vert \Delta')\to C(X\vert \Delta)$ such that $c_{(X\vert \Delta)}\circ u=c_u\circ c_{(X'\vert \Delta')}$. 

4. If $u$ is finite and orbifolde-\'etale, $c_u$ is generically finite.

5. If $(f\vert \Delta):(X\vert \Delta)\dasharrow Y$ is a fibration with general orbifold fibre special (resp. if $(g\vert \Delta):(X\vert \Delta)\dasharrow Z$ is of general type), there exists unique maps $u:Y\to C(X\vert \Delta)$ (resp. $v:C(X\vert \Delta)\dasharrow Z$ such that: $c_{(X\vert \Delta)}=u\circ f$ (resp. $g=v\circ c_{(X\vert \Delta)})$. Use lemma \ref{u}. 

The first case applies for example to the `rational quotient' $r$ of $(X\vert \Delta)$ (see \ref{rqo} below), and to its Moishezon-Iitaka fibration $M$ if $\kappa(X\vert \Delta)\geq 0$ (see \ref{tom} below).
\end{remark}

\begin{example}

1. Let $M_X:X\to \PP^1$ be the Moishezon-Iitaka fibration of example \ref{exampleiitaka}, in the notation of \S \ref{dpfa}. This is the core of $X$.

2. Let $M_X:X\to M(X)$ be the Moishezon-Iitaka fibration of any subvariety of an Abelian variety. This is the core of $X$, by \cite{U75}, thm 10.3, p.120.

3. Using the last property of remark \ref{propc} above, one sees that for any $X$, the composite map $M\circ r:X\dasharrow MR(X)$ (well-defined if $\kappa(R(X))\geq 0$, which conjecturally always holds) has special fibres, since $\Delta(r)=0$, by \cite{GHS}. Thus  $M\circ r$ factorises $c_X$ (ie: $c_X=u\circ (M\circ r)$ for some $u:MR(X)\dasharrow C(X))$. 

We will now generalise this factorisation.
\end {example}

\section{\bf The Decomposition of the Core}\label{decc}

Even if we start with some $X$ without orbifold structure (ie: $\Delta=0)$, we are, in general, immediately faced with the the existence of nonzero $\Delta's$ in the orbifold bases of the decomposition process $(M\circ r)^n$ described above (in \S \ref{dpfa}), which needs to be run in the orbifold category. So we state the decomposition structure in this later category. We need to define the maps $r$ and $M$ in this orbifold context.

We always denote by $(X\vert\Delta)$ a smooth projective orbifold.

\

The definition of $M$ for orbifolds does not present any new difficulty:

\begin{proposition}\label{tom} Assume that $\kappa(X\vert\Delta)\geq 0$. Then, there exists a (birationally) unique fibration: $M=M_{(X\vert\Delta)}:(X\vert\Delta)\to M(X\vert\Delta)$ such that:

1. Its orbifold fibres have $\kappa=0$.

2. $dim(M_{(X\vert\Delta)})=\kappa(X\vert\Delta)$.\end{proposition}

The construction is the same as when $\Delta=0$.

\begin{remark}\label{rM} We do not know if the generic orbifold fibres or orbifold base of `neat' models of $M$ are birational invariants. 
\end{remark}

\subsection{Weak Orbifold Rational Quotient.}

\

\

We shall now define an orbifold weak version of the `rational quotient' \ref{rq} in the orbifold context, assuming $C_{n,m}^{orb}$.

\begin{theorem}\label{rqo} Let $(X\vert\Delta)$ be a smooth geometric orbifold. There exists a (birationally) unique map $r=r_{(X\vert\Delta)}:X\dasharrow R=R(X\vert\Delta)$ such that:

$\bullet$ $\kappa(r\vert\Delta)\geq 0$. (assuming $C_{n,m}^{orb}$)

$\bullet$ Its generic orbifold fibres $(X_s\vert\Delta_s)$ satisfy: $\kappa_+(X_s\vert\Delta_s)=-\infty$.
(The generic orbifold fibres are well-defined since $r$ is almost holomorphic, by the first property)

The map $r$ is called the `$\kappa$-rational quotient' of $(X\vert\Delta)$.
\end{theorem}

{\bf Idea of proof:} Take a dominant connected map $f:X\dasharrow Y$ with $\dim Y$ maximal such that $\kappa(f\vert \Delta)\geq 0$. If the generic orbifold fibre has not $\kappa_+=-\infty$, one can construct a new map $g:X\dasharrow Z$ with $\dim Z>\dim Y$ such that $\kappa(g\vert \Delta)\geq 0$, using $C_{n,m}^{orb}$ $\square$

\begin{remark}\label{jr}

1. For any smooth orbifold $(X\vert\Delta)$, the composite map $M\circ r$ is always defined, since $\kappa(r\vert \Delta)\geq 0$ (the existence part assumes $C_{n,m}^{orb})$.

2. $M\circ r$ is the identity map of $X$ if and only if $(X\vert\Delta)$ is of general type.

3.\label{rr} We do not know if the orbifold base of `neat' models of $r$ are birational invariants. 

4. The map $r$ is almost holomorphic if $(X\vert \Delta)$ is finite, but not necessarily if $(X\vert \Delta)$ is logarithmic with lc, but not klt singularities. See example \ref{ernh}.
\end{remark}

\begin{example} Consider the special case where $X$ is a smooth surface, and $\Delta=0$: $r$ is the identity map if and only if $\kappa(X)\geq 0$, and has image a curve $B$ if and only if $\kappa(X)=-\infty$ and $q(X)>0$, it is the constant map if and only if $X$ is rational. 

The map $M\circ r$ is the identity map if and only if $\kappa(X)=2$, it maps to a curve if and only if either $\kappa(X)=1$, or $\kappa(X)=-\infty$ and $q(X)\geq 2$. Next, $M\circ r$ is the constant map if and only if either $\kappa(X)=0$, or $\kappa(X)=-\infty$, and $q(X)\leq 1$.

Thus $(M\circ r)=(M\circ r)^2$ is the core of $X$, except when $\kappa(X)=1$ and the orbifold base $(B\vert \Delta(M))$ of $M:X\to B$ is not of general type.  

We then have either: $\kappa(B\vert \Delta(M))=-\infty$, or: $\kappa(B\vert \Delta(M))=0$. In the first case, $r\circ M\circ r$ is the constant map. In the second case, $(M\circ r)^2$ is the constant map. In both cases $X$ is special, and $(M\circ r)^2=c$, the core of $X$. 
\end{example}

This example generalises as follows.

\subsection{The Conditional Decomposition $c=(M\circ r)^n$.}

\begin{theorem}\label{deccor} Assume $C_{n,m}^{orb}$. For any smooth projective $n$-dimensional $(X\vert\Delta)$, one has: $c=(M\circ r)^n$. 

Here $c$ is the core, $r$ the $\kappa$-rational quotient, and $M$ the orbifold Moishezon-Iitaka fibration, respectively defined in \ref{rqo} and \ref{tom}. \end{theorem}

{\bf Idea of proof:} After remark \ref{jr}, for any $k\geq 0$ and any smooth $(X\vert \Delta)$, the map $(M\circ r)^k$ can be defined. Moreover, if $(M\circ r)^{k+1}=(J\circ r)^k$ (which holds for $k=n$, by reasons of dimension), then $(M\circ r)^k$ is a fibration of general type. Thus $(M\circ r)^n$ is a fibration of general type. By the theorem \ref{extspec} below, it has also (by induction on $k$) special orbifold fibres, and the asserted equality is established (by uniqueness of $c)$ $\square$

\begin{theorem}\label{extspec} Let $(X\vert\Delta)$ be smooth, and let $f:X\to Y$ be a neat fibration. If the orbifold base and the general orbifold fibre of $f$ are special, then $(X\vert\Delta)$ is special.
\end{theorem}

\begin{remark}\label{rkextspec} This result is {\it very} false in the category of manifolds (just considering fibres and usual base, as seen in the example \ref{exampleiitaka}). It shows one of the essential improvements of working in the orbifold category. Another similar `additivity' result is the exact sequence (\cite{Ca07}, \S. 12.) of orbifold fundamental groups for fibrations.
\end{remark}

{\bf Idea of proof (of \ref{extspec}):} Assume not. Let $c:(X\vert\Delta)\to C$ be the core, assuming $dim(C)>0$. The orbifold fibres of $f$ being special are contained in those of $c$ (which are maximum for this property). Hence a factorisation $c=g\circ f$, for some $g:Y\to C$. On suitable birational models, this induces an orbifold  morphism: $g: (Y\vert \Delta_Y)\to (C\vert \Delta_{(c\vert \Delta)})$. Contradicting the specialness of $(Y\vert \Delta)$ and the hypothesis: $\kappa(C\vert \Delta_{(c\vert \Delta)})=dim(C)>0$ $\square$

We can now get (conditionally) a more concrete description of special orbifolds:

\begin{corollary}\label{spec} Assume $C_{n,m}^{orb}$. A smooth projective $n$-dimensional $(X\vert\Delta)$ is special if and only if $(M\circ r)^n$ is the constant map.\end{corollary}

\begin{remark} Thus special orbifolds are exactly the `orbifold combinations' of orbifolds with either $\kappa=0$, or $\kappa_+=-\infty$ (ie: are towers of fibrations with orbifold fibres having either $\kappa=0$, or $\kappa_+=-\infty)$. Even when $\Delta=0$, the consideration of orbifold structures is essential, as shown again by example \ref{exampleiitaka}.
\end{remark}

\section{\bf Orbifold Rational Curves}\label{orc}

The objective of these notions is to formulate the conjecture \ref{cok+} below, asserting, just as when $\Delta=0$, the equivalence between the condition $\kappa_+=-\infty$ and orbifold rational connectedness. The notion of rational curve in the orbifold context is more involved than for manifolds, and their geometry seems to be more difficult to study, too. See \cite[Chap. 5]{Ca07} for more details. We shall define the notions of uniruledness and rational connectedness for smooth geometric orbifolds $(X\vert\Delta)$. The objective being to establish for these orbifold rational rational curves the same results (\cite{MM}, \cite{M},\cite{KMM92}, \cite{KMM'}) as when $\Delta=0$.

We shall here consider mainly the `divisible' orbifold rational curves, and shall give a very brief survey of the problems arising. See again \cite{Ca07},\S 5 for more details.

\subsection{Minimal Orbifold Divisors.}

\

Let $(X\vert \Delta)$ be a {\it smooth integral} projective orbifold, and let $g:C\to (X\vert \Delta)$ be a morphism from a connected smooth projective curve $C$ such that $g(C)\subsetneq Supp(\Delta)$. Write (as usual): $\Delta:=\sum_j(1-\frac{1}{m_j}).D_j$.

There thus exists a unique smallest {\it integral} orbifold divisor $\Delta_C$ on $C$ such that $g:(C\vert \Delta_C)\to (X\vert \Delta)$ is a `divisible' orbifold morphism\footnote{One can define orbifold $\Bbb Q$-rational curves (or curves which are $\Delta^{Q}$-rational) similarly. 

In this case, $m^Q_C(a)=max_{j\in J(a)}\{ \frac{m_j}{t_{j,a}}\}$, instead.}.

Explicitly, for any point $a\in C$, the $\Delta_C$-multiplicity $m_C(a)$ of the divisor $\{a\}$ on $C$ is given by:

$\bullet$ $m_C(a)=1$ if $f(a)\notin Supp(\Delta)$.

$\bullet$ $m_C(a)=lcm_{j\in J(a)}\{\frac{m_j}{gcd({m_j,t_{j,a}})}\}$, where: $J(a):=\{j\in J$ such that: $g(a)\in D_j\}$, and $g^*(D_j)=t_{j,a}.\{a\}+\dots$, $t_{j,a}$ being the order of contact of $g_*(C)$ with $D_j$ at $g(a)$ (By convention, the only multiple of $+\infty$ is itself).

\begin{example}

1. $m_C(a)=1$ iff for each $j\in J(a)$, $m_j$ divides $t_{j,a}$.

2. If $\Delta$ is logarithmic, $\Delta_C=0$ iff $g(C)$ is disjoint from $\Delta$.
\end{example}

\subsection{Orbifold rational curves}

\begin{definition}\label{dorc} A (`divisible')\footnote{Orbifold $\Bbb Q$-rational curves are defined similarly, using the orbifold divisor $\Delta^Q_C$, instead.} rational curve on $(X\vert \Delta)$ (or: a $\Delta$-rational curve on $X$) is a map $g:C\to (X\vert \Delta)$ such that $\kappa(C\vert \Delta_C)=-\infty$. (This implies that $C\cong \Bbb P^1$, and $Supp(\Delta_C)$ has at most $3$ points with multiplicities at most $(2,3,5))$. 

When $R\subset X$ is a rational curve not contained in $Supp(\Delta)$, we say that $R$ is a $\Delta$-rational curve if so is,  in the sense above, the `normalised' inclusion $g:\hat{R}\to (X\vert \Delta)$ obtained by composing the inclusion of $R$ in $X$ with the normalisation of $R$\footnote{In \cite[\S 5]{Ca07}, maps  $g:{R}\to (X\vert \Delta)$ which are not necessarily birational to their image are considered, too. They should be important to obtain more deformations.}
\end{definition}

\begin{example}\label{Q}

1. If $R$, rational, has all its orders of contact with each $D_j$ divisible by $m_j$, then $R$ is $\Delta$-rational, with: $\Delta_C=0$.

2. If $(X\vert\Delta)$ is smooth logarithmic ($\Delta=supp(\Delta))$, the $\Delta$-rational curves are the rational curves $R$ on $X$ whose normalisation meet $\Delta$ in at most $1$ point.

For example, if $X=\Bbb P^2$, and $\Delta=D$, is a projective line with infinite multiplicity, a line $L$ (resp. an irreducible conic $C$, resp. an irreducible singular cubic $Q$) is $\Delta$-rational if and only if $L\neq D$ (resp. $C$ is tangent to $D$, resp. $Q$ is cuspidal, and tangent to $D$ in its unique singular point).

3. If $X=\Bbb P^2$, and if $\Delta$ is the union of $3$ lines in general position, with arbitrary integral multiplicities $a,b,c$, then $(\Bbb P^2\vert\Delta)$ is Fano. Indeed: $3-[(1-\frac{1}{a})+(1-\frac{1}{b})+(1-\frac{1}{c})]<0.$ There are three families of $\Delta$-rational lines covering $X$, those passing through one of the three intersection points of any two of the three given lines. We shall see in  \ref{eqg} and \ref{eqg'} that, for any finite subset $E$ of $\Bbb P^2$, there exists a curve, $\Delta$-rational irreducible, containing $E$.

4. Let $X=\Bbb P^2$, and $\Delta$ be the union of $4$ lines in general position, with multiplicities $2,2,a,b$, for arbitrary integers  $2\leq a\leq b$. Thus $(\Bbb P^2\vert\Delta)$ is Fano, since $3-[(1-\frac{1}{2})+(1-\frac{1}{2})+(1-\frac{1}{a})+(1-\frac{1}{b})]=2-[(1-\frac{1}{a})+(1-\frac{1}{b})]<0$. In this case, if $a\geq 4$, only one of the three preceding families of curves consists of $\Delta$-rational curves: those passing through the intersection of the two lines of multiplicities $a$ and $b$. If we replace the two lines of multiplicity $2$ by a conic of multiplicity $2$, we get a second family of lines which are $\Delta$-rational: the family of tangents to $C$.

5. Consider now $X=\Bbb P^2$, and $\Delta$ the union of $4$ lines in general position, with multiplicities $3,3,5,7$. Thus $(\Bbb P^2\vert\Delta)$ is Fano, since  $1/3+1/3+1/5+1/7=106/105>1$.The lines which are $\Delta^Q$-rational are finite in number, since a line meeting $Supp(\Delta)$ in at least $3$ points, gets multiplicities at least $3,3,7$, and $1/3+1/3+1/7<1$. 

It does not seem obvious to produce an explicit covering family of $\Delta$-rational curves in this case. A dimension count however shows that this might be possible (even with $\Delta'=0)$, but only for large degrees, divisible by  $3.5.7=105$. Indeed: rational planes irreducibles curves of degree $d=N. 105$ depend on  $p_N:=3(d+1)-1-3=3d-1$ parameters. Having with a line only points of contact of orders all divisible by $d'$, divisor of $d$, depends on $\frac{d}{d'}.(d'-1)= d.(1-\frac{1}{d'})$ conditions.

Thus, in our case, we get (with $d'=3,3,5,7$ successively), a number of conditions in total equal to $c_N:=d.[(1-\frac{1}{3})+(1-\frac{1}{3})+(1-\frac{1}{5})+(1-\frac{1}{7})]=d.(3-\frac{1}{105})=3.d-3.N$. We can thus expect to have a family of such curves depending on: $p_N-c_N=3N-1$ parameters of $\Delta$-rational curves (without orbifold structure, even).

See also the example \ref{ex} below.

6. The preceding example can be considered with multiplicities $(a\leq b\leq c\leq d)$ instead of $(3,3,5,7)$. The degree of the canonical bundle is then $\delta:=-3+[(1-\frac{1}{a})+(1-\frac{1}{b})+(1-\frac{1}{c})+(1-\frac{1}{d})]$, so that $(\Bbb P^2\vert \Delta)$ is Fano if and only if: $-\delta=\frac{1}{a}+\frac{1}{b}+\frac{1}{b}+\frac{1}{d}-1>0$, for example if $(a,b,c,d)=(2,3,7,41)$. Notice that, in this last case, the orbifold considered is simply-connected, and has, in particular, no orbifold \'etale cover.

The dimension count above still applies in this context and shows that rational curves of degree $d=N.m$, where $m$ is divisible by $ppcm(a,b,c,d)$, should depend on, at leat, $d.m.(-\delta)-1$ parameters. Observe that $d.(-\delta)-1=-(K_{\Bbb P^2}+\Delta).R+2-3$ if $R$ is a curve of degree $d.m$. The right-hand side of this last equality is, in fact, the Euler characteristic of the sheaf of the lift of the orbifold tangent sheaf to $R$, which computes the obstructions to deforming a orbifold morphism $\Bbb P^1\to (X\vert \Delta)$ with direct image $R$. See \cite{Ca07}, remark 5.8 for more details.

7. Let $X$ be a smooth toric projective manifold. Let $D$ be its toric  anticanonical divisor (it has normal crossings). Equip each of the components of $D$ of a finite multiplicity, and let $\Delta$ be the resulting orbifold divisor. If $R$ is a rational toric curve (closure of the orbit of a one-parameter algebraic subgroup of the torus acting on $X$), not contained in $D$, then $R$ meets $D$ in at most $2$ points in which it is unibranch. It is thus a $\Delta$-rational  curve.

8. By contrast, consider the logarithmic smooth orbifolds $(\Bbb P^2\vert \Delta),\Delta:=L+L'$ and $(\Bbb P^2\vert \Delta'), \Delta':=C$, where $L,L'$ are two distinct lines, and $C$ is a smooth conic. Both are Fano. The only $\Delta$-rational curves are, however, the lines through the intersection point $(L\cap L')$, while there exists, for every $d>0$, an explicit $(d+1)$-dimensional family of $\Delta'$-rational curves of degree $d$. The reason for this different behaviour lies in the lc, but not klt singularity of $(\Bbb P^2\vert \Delta)$. See \cite{Ca07}, example  5.14.

\end{example}

\subsection{Orbifold Uniruledness and Rational Connectedness}

\
  \begin{definition} A smooth $(X\vert \Delta)$ is {\bf uniruled (resp. R.C)} if any generic $x$ (resp. $(x,y))$ in $X$ is contained in a $\Delta$-rational curve.
  \end{definition}

  \begin{example}

1. Let $(X\vert\Delta)$ be the toric example \ref{Q}.7 above. Then $(X\vert \Delta)$ is rationally connected, since all of its toric rational curves (those not contained in $D$) are $\Delta$-rational.

2. The orbifolds $(\Bbb P_2\vert \Delta)$ of examples \ref{Q}.(3+4)  are $RC$, by the explicit covering families of $\Delta$-rational curves displayed there.

3. The orbifolds $(\Bbb P_2\vert \Delta)$ of examples \ref{Q}.5, 6 should be $RC$, by  the counting argument given there.

\end{example}

   Some properties of rationally connected manifolds extend immediately:
   
   \begin{proposition} Assume $(X\vert \Delta)$ is rationally connected. Then:
   
   1. $\pi_1(X\vert \Delta)$ is finite.
   
   2. $H^0(X,S^N(\Omega^p(X\vert \Delta)))=0$, for any $N,p>0$. Thus $\kappa_{++}(X\vert \Delta)=-\infty$.
   \end{proposition}

 See \cite{Ca07} for a proof.
 
 \begin{remark} It is not obvious (if true) that orbifold uniruledness and rational connectedness are birationally invariant properties (for finite, integral, smooth) orbifolds. Indeed, some orbifold rational curves may not lift under an (orbifold) blow-up. An example is given in \cite[5.8]{Ca07}.
 \end{remark}

 \subsection{Fano Orbifolds}
 
 \

 The following question is the easiest one of a series of questions asking whether results known when $\Delta=0$ extend to the orbifold context:
 
 \begin{question} Let $(X\vert \Delta)$ be a smooth, integral, finite geometric orbifold, with $X$ projective. Assume $(X\vert \Delta)$ is Fano (ie: $-(K_X+\Delta)$ is ample). Is then $(X\vert \Delta)$ uniruled (resp. rationally connected)? Is this at least true when $Pic(X)\cong \Bbb Z$?
  \end{question}

  \
  
  \begin{example} Assume that $(X\vert \Delta)$ is smooth, finite and Fano, with $X$ projective and $Pic(X)\cong \Bbb Z$ generated by the ample line bundle $H$. Let $r\leq (n+1)$ be the divisibility index\footnote{That is, the largest positive integer $s$ such that $K_X\in \{s.Pic(X)\}$.} of $X$. Assume that $\Delta=(1-\frac{1}{m}).D$, with $D$ a member of $\vert k.H\vert$ with only normal crossings as singularities. Thus $(1-\frac{1}{m}).k<r$, since $(X\vert \Delta)$ is supposed to be Fano. Assume, more strongly, that $k\leq r$. It follows then from \cite[Theorem 5. 8, p. 28]{KM} that there exists a covering family of rational curves on $X$ meeting $D$ in at most two points (after normalisation). In other words: $(X\vert D)$ is uniruled. Thus so is $(X\vert \Delta)$. The assumption that $k\leq r$ is met in particular when $r=1$. The situation thus seems to be more involved when $r$ increases, in particular when it is maximum, equal to $(n+1)$, that is when $X\cong \Bbb P^n$.
  \end{example}

  \
  
  Other questions ask whether Miyaoka-Mori's Bend and Break, Miyaoka's generic semi-positivity, Graber-Harris-Starr extend to the orbifold context as well. See \cite[\S5]{Ca07} for details, and corollary \ref{comp} for a very partial positive answer.

\subsection{Orbifold Uniruledness and Canonical Dimension}

\

\

A fundamental conjecture is:

 \begin{conjecture}\label{cok+} 1. $(X\vert \Delta)$ is uniruled if $\kappa(X\vert \Delta)=-\infty$
 
 2. $(X\vert \Delta)$ is rationally connected if $\kappa_+(X\vert \Delta)=-\infty$.
 \end{conjecture}

    The converse statements are easily shown to be true. Contrary to the case $\Delta=0$, it is not known that the part 1 of this conjecture implies its part 2.
   
   The only known case is when $n=2$ and $\Delta$ logarithmic, by \cite{KM98}. It were interesting (and probably feasible) to show that their result holds for finite multiplicities as well.

\begin{example}\label{ex} Let $(\Bbb P^2\vert \Delta)$ be the Fano examples of \ref{Q}(5). It should be $RC$, by the preceding conjecture. See \ref{Q}(5) for a possible direct verification.It were interesting to have a deformation-theoretic proof of this property. Notice that such parametrised orbifold rational curves would give, for generic complex numbers $u,v,w$, a solution to the equation $u.P^3+v.Q^3+w.R^5=S^7$ for complex polynomials $P,Q,R,S$ of respective degrees at most $35,35,21,15$. \end{example}

\subsection{Global Quotients: Lifting and Images of Rational Curves}

 \begin{definition} Let $(X\vert \Delta)$ be smooth. Let $g:X'\to (X\vert \Delta)$ be finite,  with $X'$ smooth, and ramifying {\bf at least} (resp. {\bf exactly}) above $\Delta$. 
 
 This means that $g$ is \'etale over the complement of $Supp(\Delta)$, and ramifies at order $m'_j$, with $m'_j$ a multiple of $m_j$ (resp. $m'_j=m_j)$ above each $D_j$.
 
 We say that $(X\vert \Delta)$ is a {\bf global quotient} if  if there exists $g:X'\to (X\vert \Delta)$ with $X'$ smooth, $g$ finite and ramifying exactly above $\Delta$.
 \end{definition}

 \begin{example}\label{eqg} Let $X=\Bbb P^2$, $\Delta=\sum_{j=1}^{j=3}(1-\frac{1}{m_j}).D_j$, with integral $m_j>1$, and $D_j$ the lines of equation $T_j=0$, in the  homogeneous coordinates $(T_1,T_2,T_3)$. Set $m:=lmc\{m_1,m_2,m_3\}$. Let $f:\Bbb P^2\to \Bbb P^2$ be defined by $f(U_j)=T_j:=U_j^m$, for $j=1,2,3$. Then this morphism ramifies at least (resp.exactly) above $\Delta$ for any choice of the $m_j's$ (resp. if $m_j=m, \forall j)$.
\end{example}

\begin{remark} If $(X\vert \Delta)$ is a global quotient with $\Delta\neq 0$, then it is not simply-connected (more precisely, its orbifold fundamental group must have order divisible by $ppcm(m_j's))$. For example, $\Bbb P^2$ with $\Delta$ consisting of $4$ lines with multiplicities either $(3,3,5,7)$, or $(2,3,7,41)$ is not a global quotient: the orbifold fundamental groups have cardinalities $3$ and $1$ respectively.
\end{remark}

  \begin{proposition} (\cite{Ca07}) Assume that $(X\vert \Delta)$ is a global quotient, as above. Then:
  
  1. For each rational curve $R'\subset X'$, $R:=g(R')$ is $\Delta$-rational.
  
  2. For each $\Delta$-rational $R\subset X$, each component $R'$ of $g^{-1}(R)$ is rational. 
     
  \end{proposition}
    
\begin{corollary}\label{comp} If $(X\vert \Delta)$ is a global quotient, we have (among many other similar statements):
   
   1. $(X\vert \Delta)$ is uniruled (resp. R.C) iff so is $X'$. 
   
   2. $(X\vert \Delta)$ is RC iff it is RCC (i.e: if any two generic points of $X$ are connected by a chain of $\Delta$-rational curves).
   
   3. If $(X\vert \Delta)$ is Fano, it is RC.
   
   4. If $C\subset X$ is a curve s.t: $-(K_X+\Delta).C<0$ through $a\in (X-Supp(\Delta))$, there exists  a $\Delta$-rational curve through $a$.

   \end{corollary}
   
   \begin{example}\label{eqg'} Let $(X\vert\Delta)$ be the example \ref{eqg}. By the above \ref{comp}, $(X\vert\Delta)$ is $RC$. Notice that such $\Delta$-rational curves are usually highly singular.
\end{example}

       \begin{remark} The situation above is very rare. Its interest is however to show that if the deformation theory of rational curves on smooth DM stacks could be developed to the same point as on manifolds, the same results could be obtained for $\Delta$-rational curves on smooth orbifolds, not necessarily global quotients, as when $\Delta=0$, by replacing the manifold $X'$ above by the smooth DM stack $\sX\to X$ associated to $(X\vert \Delta)$. The following example shows however that some new ideas or techniques seem to be needed.  \end{remark}

\begin{example} Miyaoka-Mori's Bend and Break does, however, not hold in general for arbitrary smooth, finite, integral projective orbifolds $(X\vert \Delta)$. Let us indeed choose an isotrivial family of elliptic smooth plane cubics going through a point $a\in \Bbb P^2$, and degenerating to a union of $3$ lines $D'_j,j=1,2,3$ (take for example $x^3+y^3+s=0$ in affine coordinates $(x,y), s$ being a parameter). Let us now blow-up the line $D'_1$, containing the point $a$, in $3$ (or more) of its generic points, none of which is $a$, obtaining so the manifold $X$, with three (or more) $(-1)$-curves $E_k,k=1,2,3$. We equip now each of the $E_k's$ with a large multiplicity $m$ (at least $3$ is sufficient), obtaining so the orbifold divisor $\Delta:=(1-\frac{1}{m}). (E_1+E_2+E_3)$ on $X$. The strict transform $D_1$ in $X$ of the line $D'_1$ is thus the only rational curve produced by Miyaoka-Mori's Bend and Break process through the point $a$. But this is not a $\Delta$-rational curve, since it intersects transversally each of the $E_k's$. Observe that the orbifold $(X\vert \Delta)$ is not a global quotient, since the complement of the support of $\Delta$ is simply-connected.
\end{example}

\subsection{$\Delta$-rational curves vs $\Delta^Q$-rational curves.}

\

Let $(X\vert \Delta)$ be a smooth {\it integral} projective orbifold. We have two sets of rational curves on $X$: first the set $Ratl^Z(X\vert \Delta)$ consisting of the `divisible' ones, and the set (see footnote in \S 9.A for the definition)  $Ratl^Q(X\vert \Delta)$ of `$\Delta^Q$-rational curves' on $X$ (with respect to which one can define as in \label{durc} the notions of $Q$-uniruledness and $Q$-rational connectedness). Obviously, $Ratl^Z(X\vert \Delta)\subset Ratl^Q(X\vert \Delta)$. Exemples easily show that the inclusion is strict in general. However, we have the following: 

\begin{conjecture}\label{zvsqrat} 

1. If $(X\vert \Delta)$ is $Q$-uniruled (resp. $Q$-RC) , it is uniruled (resp. $RC)$.

2. If $(X\vert \Delta)$ is of general type, there exists a closed proper algebraic subset of $X$ containing all $Q$-rational curves of $(X\vert \Delta)$.\end{conjecture}

Concerning \ref{zvsqrat}.(2), see \cite{PR} for the case of general hypersurfaces of the projective space, and \cite{FB} for the surface case, by an arithmetic approach. (I thank A. Levin for this reference and an interesting discussion on this topic, leading to the second part of the conjecture \ref{zvsqrat} above).

\

\section{\bf Some relationships to LMMP and Abundance Conjecture.}\label{abund}

\

The objective is to show how to deduce (when $X$ is projective) the conjectures $C_{n,m}^{orb}$ and \ref{cok++} respectively from the two standard conjectures of the LMMP: a weak version of the abundance conjecture, and the existence of log-minimal models in the log-canonical case. A third property is actually needed: the stability of the logarithmic cotangent sheaf for log-canonical pairs with first Chern class zero. I thank M. P\u aun for explaining me the notion of numerical dimension.

\

We assume in this section that $X$ is a complex projective connected, $n$-dimensional $\Bbb Q$-factorial normal space, $A$ and $D$ be a $\Bbb Q$-divisors on $X$, with $A$ ample. 

The {\it numerical dimension of $D$} is defined as being the following invariant: $\nu(X,D):=limsup_{m>0}\frac{Log(h^0(X,mD+A))}{Log(m)},$ for $m>0$ integral and sufficiently divisible. 

Easy standards arguments show that:

1. $\nu(X,D)\in \{-\infty, 0,1,\dots,n\}$.

2. $\nu(X,D)$ does not depend on the choice of $A$. 

3. $\nu(X,D)\geq \kappa(X,D)$.

4. $\nu(X,D)=-\infty$ if and only $D$ is not pseudo-effective (this is one of the definitions of pseudo-effectivity).

5. $\nu(X,D)=\kappa(X,D)$ if $\kappa(X,D)\geq 0$ (ie: if $D$ is effective).

\

One (weak) form of the so-called `Abundance Conjecture' is the following:

\begin{conjecture}\label{abcj} Assume $(X\vert D)$ is a `log-canonical pair'. Then:

 $\nu(X,K_X+\Delta)=\kappa(X,K_X+\Delta)$.
\end{conjecture}

This is known\footnote{Even for $\Bbb R$-divisors.} when $D$ is `big' and $(X\vert D)$ is klt (\cite{BCHM},\cite{Pa}). This is also known when $\nu(X,K_X+\Delta)=0$ if $q(X)=0$, it follows from \cite{N}and \cite{Bo}; when moreover $\Delta=0$, the general case follows from \cite[\S 3]{CP05}.

\begin{proposition} Assume conjecture \ref{abcj}. Then Conjecture $C_{n,m}^{orb}$ is true.
\end{proposition}

{\bf Proof:} Let $f:X\to Y$ be a neat fibration with $(X\vert \Delta)$ a smooth orbifold, with $\kappa(X_y,K_{X_y}+\Delta_y)\geq 0$, for $y\in Y$ general. Let $(Y\vert \Delta_Y)$ be the orbifold base of $(f\vert \Delta)$ (ie: $\Delta_{(f\vert \Delta)}:=\Delta_Y)$. Then $f_*(m(K_{X/(Y\vert \Delta_Y)}+\Delta))$ is weakly-positive for $m$ large and divisible, by theorem \ref{vieh'}. We assume that $\kappa(Y,K_Y+\Delta_{(f\vert \Delta)})\geq 0$, otherwise the statement is trivial. 

Thus $m(K_{X}+\Delta)+f^*(B)=m(K_{X/(Y\vert \Delta_Y)}+\Delta)+f^*(K_Y+\Delta_Y)+f^*(B)$ is effective for any $\Bbb Q$-ample divisor $B$ on $Y$. And so $K_X+\Delta$ is pseudo-effective, hence effective, by conjecture \ref{abcj}. It then follows from easy arguments (see for example \cite[2.3-4,p. 516]{AC}) that $C_{n,m}^{orb}$ subadditivity holds in this case $\square$

\

The standard second conjecture of the LMMP is (see \cite[conjecture 1.1]{Bi}, for example):

\begin{conjecture}\label{lmmp}Let $(X\vert\Delta)$ be a log-canonical pair. There exists a sequence of divisorial contractions and flips $s:X\dasharrow X'$ such that, if $s_*(\Delta):=\Delta'$, then $(X'\vert\Delta')$ is log-canonical, $\Bbb Q$-factorial, and either:

1.  $K_{X'}+\Delta'$ is nef, or:

2. There exists a fibration $f:X'\to Y'$ with $-(K_{X'}+\Delta')$ relatively ample.
\end{conjecture}

This conjecture is known in the klt case, by \cite{BCHM}, and also when $n:=dim(X)\leq 3$, by \cite{Sh}, see also \cite{Bi} and the references there.

\begin{conjecture}\label{c1=0} Assume $(X\vert \Delta)$ is a `log-canonical pair' with $c_1(X\vert \Delta)=0$. Then: 

1. $\kappa_{++}(X\vert \Delta)=0$, defining the left-hand side as being $\kappa_{++}(X'\vert \Delta')$, for any log-resolution $(X'\vert \Delta')$ of $(X\vert \Delta)$. 
This is implied by the following orbifold version of Miyaoka's semipositivity theorem: 

2. if $C$ is the complete intersection of large multiples of given ample divisors $H_i$ of $X$, and contained in the open set where $(X\vert \Delta)$ is smooth, then the restriction of $S^N\Omega^p(X\vert \Delta)$ to $C$ should be $(H_1,H_2,\dots,H_{n-1})$-semi-stable, for any $p>0$ and $N>0$.
\end{conjecture}

This conjecture is known when $\Delta=0$, using either the existence of Ricci-Flat metrics (constructed by S.T. Yau), or Miyaoka's generic semi-positivity theorem. The same first approach seems accessible in the klt case. The log-canonical case might require new ideas. The orbifold version of Miyaoka semipositivity theorem might be shown by extending the arguments to the orbifold case, using the orbifold rational curves introduced in \S\ref{orc}.

\begin{theorem} Assume that conjectures $C_{n,m}^{orb}$, \ref{c1=0}, and \ref{lmmp} hold. Then conjecture \ref{cok++} holds. 

Thus conjecture \ref{cok++} follows from conjectures \ref{abcj}, \ref{c1=0}, and \ref{lmmp} of the LMMP.

When $(X\vert \Delta)$ is smooth, finite with $X$ projective, conjecture \ref{cok++} follows from conjectures \ref{abcj} and \ref{c1=0}, since conjecture \ref{lmmp} is known from \cite{BCHM}.
\end{theorem}

{\bf Proof:} Let $(X\vert \Delta)$ be a smooth orbifold with $X$ projective. Combining conjectures \ref{c1=0} and \ref{lmmp}, we see that conjecture \ref{cok++} holds when $\kappa(X\vert \Delta)=0$. Let $\sF\subset \Omega^p_X$ be a rank-one coherent subsheaf.

There are now $2$ exclusive cases: either $\kappa(X\vert\Delta)\geq 0$, or $\kappa(X\vert\Delta)=-\infty$.

In the first case, let $f:X\to Y$ be the Iitaka fibration associated to $(K_X+\Delta)$. We can assume that $f$ is regular, by making an orbifold elementary modification, which does not change the birational invariants we are interested in. Thus $dim(Y)=\kappa(X\vert \Delta)$, and $\kappa(X_y\vert \Delta_y)=0$ over the generic fibre $X_y$ of $f$. By a further orbifold modification, we can assume that the saturation of $\sF^{\otimes N}$ in $S^N\Omega^p(X\vert \Delta)$ is locally free for a suitable $N$ such that $\kappa(X\vert \Delta,\sF)$ is given by the linear system associated to this saturation. Since the normal bundle to $X_y$ in $X$ is trivial, and since conjecture \ref{cok++} holds for $(X_y\vert \Delta_y)$, by the assumptions made, we conclude that $\kappa(X_y,\sF)\leq 0$. By considering $f_*(\sF^{\otimes N})$, we see that $\kappa(X\vert \Delta,\sF)\leq dim(Y)=\kappa(X\vert \Delta)$. 

It remains to deal with the case when $\kappa(X\vert \Delta)=-\infty$. Applying conjecture \ref{lmmp}, we get a birational sequence $s:(X\vert \Delta)\dasharrow (X'\vert \Delta')$ and a log-Fano fibration $f:(X'\vert \Delta')\to Y$, with $dim(Y)<n$ and $(X'\vert \Delta')$ log-canonical. Let $H'$ be a general member of $m.(-(K_{X'/Y}+\Delta')+f^*(A))$, for $A$ sufficiently ample on $Y$, and $m$ sufficiently large, so that $(X'\vert \Delta'+\frac{1}{m}.H'):=(X'\vert \Delta")$ is log-canonical, with $(K_{X'_y}+\Delta"_y)$ trivial on the generic fibre $X_y$ of $f$. Choosing $m$ sufficiently large, and making a suitable orbifold modification of $(X\vert \Delta)$, we can moreover assume that $s:(X\vert \Delta)\to (X'\vert \Delta')$ and $s:(X\vert \Delta^+)\to (X'\vert \Delta")$ are regular and log-resolutions, where $\Delta^+:=\Delta+\frac{1}{m}.H$, and $H$ being the strict transform of $H'$ in $X$. 

Considering the restriction $s_y:X_y\to X'_y$ of $s$ over a generic $y\in Y$, we see that $s_y:(X_y\vert \Delta_y)\to (X'_y\vert \Delta'_y)$ is a log-resolution with $-(K_{X'_y}+ \Delta'_y)$ ample on $X'_y$. From Lemma \ref{fanok++} below we deduce that $\kappa_{++}(X_y\vert \Delta_y)=-\infty$ , and lemma \ref{desc} then shows that any non-zero section $w$ of $\sF$ is of the form $w=f^*(v)$, for some section $v$ of $S^N\Omega^q(Y'\vert \Delta_{Y'})$, where $(Y'\vert \Delta_{Y'})$ is the smooth orbifold base of some suitable neat model of $f\circ s:(X\vert \Delta)\to Y$, and $0\leq q\leq p$. Finally, lemma \ref{fanofib} below implies that $q=0$, that is: $\kappa_{++}(X\vert \Delta)=-\infty$ if $\kappa_+(X\vert \Delta)=-\infty$, which is the claimed implication $\square$

\begin{lemma}\label{fanok++} Let $g:(X\vert \Delta)\to (X'\vert \Delta')$ be a birational map from the smooth orbifold $(X\vert \Delta)$ to the log-canonical Fano pair $(X'\vert \Delta')$ such that $f_*(\Delta)=\Delta'$. Assume that conjecture \ref{c1=0}.(2) holds. Then, for any polarisations of $X'$, and any general Mehta-Ramanathan curve $C\subset X'$, identified with its strict transform on $X$, the following properties hold:

1. For any finite sequence of pairs of positive integers $(N_h,q_h),h=1,\dots, t$, and any subsheaf $\sF\subset S^{N_1}\Omega^{q_1}(X\vert \Delta)\otimes \dots \otimes S^{N_t}\Omega^{q_t}(X\vert \Delta)$, the restriction $\sF_{\vert C}$ has negative degree. In particular, $H^0(X,S^{N_1}\Omega^{q_1}(X\vert \Delta)\otimes \dots \otimes S^{N_t}\Omega^{q_t}(X\vert \Delta))=\{0\}$.

2. $\kappa_{++}(X\vert \Delta)=0$.
\end{lemma}

{\bf Proof:} Let $\sG:=S^{N_1}\Omega^{q_1}(X\vert \Delta)\otimes \dots \otimes S^{N_t}\Omega^{q_t}(X\vert \Delta)$, and let $L$ be a line bundle of degree $0$ on $C$. It is sufficient (by considering $\Lambda ^r\sG$, for any $r>0)$ to show that the degree of any rank one coherent subsheaf of $\sG$ is negative over $C$, and even that $H^0(C,\sG\otimes L)=0$. Let $H'=\sum_{j=1}^{j=n}H_j$, where the $H_j's$ are general members of $m.(-(K_{X'}+\Delta'))$, $m>0$ being a sufficiently large integer, the $H_j's$ being chosen so that $(X'\vert \Delta'+\frac{1}{mn}.H'):=(X'\vert \Delta")$ is log-canonical, with $(K_{X'}+\Delta")$ trivial on $X'$, and such that $\Delta^+:=(\Delta+\frac{1}{mn}.H)$ has normal crossings support\footnote{Thus $\Delta^+$ is not integral (unless $m=n=1)$, even if $\Delta$ is. So $\Q$-multiplicities (here equal to $\mu=\frac{mn}{mn-1}=1+\frac{1}{mn-1})$ are needed in conjecture \ref{c1=0} in the present argument.}, $H$ being the strict transform of $H'$ in $X$. Choose $C\subset X'$ to be a generic curve cut out by Mehta-Ramanathan very ample linear systems, and meeting $H'$ transversally, but not meeting the indeterminacy locus of $g^{-1}$,and so identified with its strict transform on $X$, then any non-zero section of $L\otimes\sG^+_{C}:=L\otimes S^{N_1}\Omega^{q_1}(X\vert \Delta^+)\otimes \dots \otimes S^{N_t}\Omega^{q_t}(X\vert \Delta^+)_{\vert C}$ has no zero at all on $C$, since it generates a rank one subsheaf of the locally free sheaf $L\otimes\sG^+_{C}$, assumed to be semi-stable, and with trivial determinant.

Assume now that $w$ is a non-zero section of $L\otimes \sG$. Replacing $w$ by a suitable tensor power, we may assume that $N:=\sum_{j=1}^{j=t}N_h.q_h\geq \frac{n}{p.(1-\frac{1}{mn})}$. Because of the natural inclusion $\sG\subset \sG^+$, $w$ is also a section of $L\otimes\sG^+_{\vert C}$, and has thus no zero (as a section of $L\otimes\sG^+_{\vert C})$. Now choose  a generic point $a\in X'$ such that $w(a)\neq 0$ as a section of $L\otimes\sG_{\vert C}$. We choose $a$ outside of the support of $\Delta'$ and belonging to the smooth locus of $X'$. On the other hand, since the $H_j$ can be chosen arbitrary generically, we may assume that they build a system of coordinate hyperplanes for suitable local coordinates at $a$. Now we check immediately, using lemma \ref{diffsheaf} and multilinearity, that in the local $\sO$$_{X'}$-basis $dz_{(J)}$ of $\sG^+$ described in \S\ref{inva}, all coordinates of $w$ vanish at $a$. This contradicts the non-vanishing property above, and proves the claims $\square$

\begin{lemma}\label{desc} Let $(X\vert \Delta)$ be a smooth Fano orbifold, and $f:X\to Y$ be a fibration with generic smooth orbifold fibre $(X_y\vert \Delta_y)$. 

Assume that, for any finite sequence of pairs of positive integers $(N_h,q_h),h=1,\dots, t$, one has: $H^0(X_y,S^{N_1}\Omega^{q_1}(X_y\vert \Delta_y)\otimes \dots \otimes S^{N_t}\Omega^{q_t}(X_y\vert \Delta_y))=\{0\}$.

Then any section $w$ of $S^N\Omega^p(X\vert \Delta)$ is of the form $w=f^*(v)$, for some section $v$ of $S^N\Omega^q(Y'\vert \Delta_{Y'})$, where $(Y'\vert \Delta_{Y'})$ is the smooth orbifold base of some suitable neat model of $f\circ s:(X\vert \Delta)\to Y$, and for some $q$ with $0\leq q\leq p$.

In particular: $\kappa_{++}(X\vert \Delta)\leq dim(Y)$.
\end{lemma}

{\bf Proof:}  This is just Lemma \ref{f*surj} $\square$

\begin{lemma}\label{fanofib} Let $(X\vert \Delta)$ be a smooth orbifold with $\kappa_+(X\vert \Delta)=-\infty$, and let $f:X\to Y$ be a fibration. If $(Y'\vert \Delta_{Y'})$ be the orbifold base of some suitable neat model of $f$ such that $(Y'\vert \Delta_{Y'})$ is smooth, then $\kappa_+(Y'\vert \Delta_{Y'})=-\infty$.
\end{lemma}

{\bf Proof:} After an orbifold modification of $(X\vert \Delta)$, we can assume that $f$ is an orbifold morphism, in which case $f^*(S^N\Omega^p(Y'\vert \Delta_{Y'}))\to S^N\Omega^p(X\vert \Delta)$ is well-defined for any $N,p\geq 0$, so that: $\kappa_+(Y'\vert \Delta_{Y'})\leq \kappa_+(X\vert \Delta)=-\infty$ $\square$

\begin{corollary}\label{knownck++} The conjectures \ref{abcj}, \ref{lmmp} and \ref{c1=0}, and thus $C_{n,m}^{orb}$ and \ref{cok++} are known (see, for example, \cite{Kal} and \cite{Bi}) in the following special cases:

1. $n\leq 2$

2. $n=3$ and $\Delta=0$
\end{corollary}

\section{\bf Some Conjectures}

\subsection{Lifting of properties}\label{conj}

\

We assume in this \S \ref{conj} the existence of the $\kappa$-rational quotient' for smooth projective orbifolds (this is true if so is $C_{n,m}^{orb})$. Thus we have the decomposition $c=(M\circ r)^n$ for any such orbifold. Then:

\begin{proposition}\label{lp} Let $P$ be a property of smooth projective orbifolds which is:

1. birational.

2. satisfied when $\kappa=0$ and $\kappa_+=-\infty$.

3.  stable by exensions (ie: satisfied by $(X\vert \Delta)$ if satisfied by the general orbifold fibre 
and the orbifold base of some neat fibration $(f\vert \Delta):(X\vert \Delta)\to Y)$.

Then $P$ is satisfied by any special orbifold.

Assume moreover that the property $P$:

4. is {\it not} satified by the positive-dimensional orbifolds of general type.

5. is satisfied by $(X\vert \Delta)$ if satisfied by $(X'\vert \Delta')$, and there exists a surjective orbifold morphism $(X'\vert \Delta')\to (X\vert \Delta)$.

Then $P$ is satisfied {\it exactly} by the special orbifolds.
\end{proposition}

\subsection{The case of special orbifolds}

\

\

They should qualitatively behave like rational and elliptic curves. They are, defined similarly: coherent rank $1$ subsheaves $L\subset \Omega^p_X$ should have $\kappa(X\vert \Delta,L)<p$, for any $p>0$. We formulate the conjectures in general, although the invariants below have not been defined when $\Delta\neq 0$ here.

\begin{conjecture}\label{cspec} Let $(X\vert \Delta)$ be a smooth projective orbifold.

$\bullet$ If $(X\vert \Delta)$ is special, then $\pi_1(X\vert \Delta)$ is almost abelian.

$\bullet$ $(X\vert \Delta)$ is special if and only if its Kobayashi pseudometric $d_{(X\vert \Delta)}\equiv 0$. 

$\bullet$ $(X\vert \Delta)$ (defined over a number field) is special if and only if potentially dense. 

$\bullet$ if $(X\vert \Delta)$ is special, so are its deformations (see \cite{Ca07} for a precise statement).

\end{conjecture}

\

\begin{remark}\label{cs}

1. The motivation for these conjectures is of course \ref{lp}, which might give a possible approach for a proof:

2. The conjecture should hold for the crucial cases: $\kappa(X\vert \Delta)=0$ and $\kappa_+(X\vert \Delta)=-\infty$, which might be proved first by other specific methods for the particular cases $K_X+\Delta<0$ and $c_1(K_X+\Delta)=0$ (In the case of $\pi_1$, for example, one may think of $L^2$ theory, and Ricci-flat metrics). See also the remark \ref{rfpm}.(1).

3. The properties should be shown to be preserved by `extensions' (ie: should hold for $(X\vert \Delta)$ if they hold for the orbifold fibres and base of a fibration $f:X\to Y)$. Notice that the consideration of the orbifold structure on the base precisely means that local obstructions to lifting vanish (in codimension one at least). The expected `extension' property is thus a kind of `local-to-global' principle. This extension property holds for the orbifold fundamental group (see \cite{Ca07},\S 12).

4. Assuming $C_{n,m}^{orb}$, the conclusion would follow from theorem \ref{deccor}.

5.  One might even wonder, in case of the second conjecture above, whether being special is not equivalent to have any two points joined by an entire (transcendental) orbifold curve $h:\C\to X$.

6. The hyperbolicity and arithmetic conjectures above have an obvious function field analogue. See \cite{Ca07} for details.
\end{remark}

\subsection{The general case}

\

The general case should be `split' into its two antithetical parts (special and `general type') by the core. See \cite{Ca07} for definitions and details.

\begin{conjecture}\label{cgen} Let $(X\vert \Delta)$ be a smooth projective orbifold, and let $c:(X\vert \Delta)\to C(X\vert \Delta)=C$ be its core. Then:

$\bullet$ $d_{(X\vert \Delta)}=(c)^*(\delta)$, with $\delta:=d_{(C\vert \Delta(c\vert \Delta))}$, the Kobayashi pseudometric on the base orbifold of $(c\vert \Delta)$, which is a metric on some Zariski dense open subset $U$ of $C$.

$\bullet$ If $(X\vert \Delta)$ is defined over a number field, then $(U\cap c((X\vert \Delta)(k)))$ is finite, for any number field of definition of $(X\vert \Delta)$, and any `model' of $(X\vert \Delta)$. 

$\bullet$ The fibration $c$ deforms under deformations of $(X\vert \Delta)$, in particular $\dim C(X\vert \Delta)$ is constant in such a deformation.

\end{conjecture}

\begin{remark} Let us motivate the first conjecture above (the second one is similar): since the general orbifold fibres of $c$ are special, their Kobayashi pseudometric vanishes, after \ref{cspec}. Thus $d_{(X\vert \Delta)}=(c_X)^*(\delta)$, for some $\delta$ on $C$. By the definition of Kobayashi pseudometric on the orbifold base of $c$, we have: $\delta\geq d_{(C\vert \Delta(c\vert \Delta))}$. The reversed inequality should then follow from the fact that, locally on the complement of a $2$-codimensional subset of $C(X)$, the orbifold maps $h:\D\to (C(X)\vert \Delta(X\vert \Delta))$ lift to $X$.

The assertion that $\delta$ is a metric generically on $C(X\vert \Delta)$ is simply an orbifold version of Lang's conjecture. Similarly for the the finiteness in the arithmetic conjecture.

\end{remark}

\subsection{Families of Canonically Polarised Manifolds.}\label{fcpm}

\

\

When $(X\vert \Delta)$ is a smooth logarithmic orbifold, the core and the notion of specialness produce new invariants of the quasi-projective manifold $U=X-\Delta$, independent from its smooth projective compactifications. This framework appears to be suitable also in moduli problems.

\begin{conjecture}\label{isotcj} {\bf (``Isotriviality Conjecture")} We conjecture (in \cite{Ca07}, \S 12.6) that algebraic families of canonically polarised manifolds parametrised by {\it special} quasi-projective varieties $X^0$ are isotrivial, and so that, for any quasi-projective base,  the moduli map factors through the `core' of the base of this base . \end{conjecture}

This extends and unifies previous conjectures by Viehweg-Zuo (\cite{VZ}) and Kebekus-Kov\`acs (see \cite{JK} and \cite[\S 12.6]{Ca07} for details). This conjecture is proved in dimensions $3$ or less in \cite{JK} (see also relevant references there).

 In fact, this conjecture can be reduced to other conjectures of classification, using the Viehweg-Zuo sheaves in \cite{VZ}:

\begin{proposition}\label{pfcm} The conjecture above (ie: algebraic families of canonically polarised manifolds parametrised by {\it special} quasi-projective varieties are isotrivial) is true when $dim(X^0)=d$ if conjectures $C_{n,m}^{orb}$ and \ref{cok++} hold for dimensions at most $(d-1)$. This conjecture is thus, by the observations made in \S\ref{abund}, a consequence in dimension $d$ of conjectures \ref{abcj},\ref{lmmp} and \ref{c1=0} in dimension at most $(d-1)$.
\end{proposition}

{\bf Proof:} It essentially consists in checking the properties 1-3 listed in proposition \ref{lp}. Let $X^0=(X-D)$ be a smooth logarithmic compactification of the base space $X^0$ of this base space. Then $c=(M\circ r)^n$ is the constant map to a point, since $X^0$ is special and $n$-dimensional. It is thus sufficient to show, inductively on $n$, that if $\mu:X^0\to \sM$ is the moduli map to the coarse moduli stack induced by the given family (and constructed in \cite{V95}, theorem 1.11), that $\mu$ factorises through $f$ if $(f\vert D):(X\vert D)\to Y$ is a fibration with general orbifold (logarithmic) fibres having either $\kappa=0$, or $\kappa_+=-\infty$. In other words: it is sufficient to show the conjecture when $\kappa(X\vert D)=0$, and when $\kappa_+(X\vert D)=-\infty$. 

We now assume one of these properties to hold, and also that  $Var(g^0)>0$, where $Var(g^0)$ is the generic rank of the map $\mu$, or equivalently, of the Kodaira-Spencer map associated to the family of canonically polarized manifolds under consideration. The main theorem of \cite{VZ} asserts the existence of a line bundle $A\subset Sym^N(\Omega^1_X(log D))$ with $\kappa((X\vert D),A)\geq Var(g^0)$. An important refinement of this result (\cite{JK}, theorem 9.3), shows that $A\subset Sym^N(B)$ for some subsheaf  $B\subset \Omega_X^1(log D)$, generically of the form $B=Im(d\mu:\mu^*(\Omega^1_{\sM})\to \Omega^1_{X^0})$. Using now conjecture \ref{cok++}, we get: $\kappa(X\vert D,A)\leq \kappa_{++}(X\vert D)=\kappa(X\vert D)\leq 0$, contradicting our hypothesis that $Var(g^0)>0$ $\square$

\begin{remark}\label{rfpm}

1. The conjecture \ref{cok++} is needed only for the cases where $\kappa_+=-\infty$ or $\kappa=0$. The first case were a consequence of conjecture \ref{cok+}, and the second might be solved using LMMP to reduce to the case of a trivial first Chern class, possibly accessible by constructing orbifold Ricci-flat K\" ahler metrics. See \S\ref{abund} for details.

2. The meaning of the above conjecture is also that any subvariety of the moduli stack $\sM$ should be of logarithmic general type. In particular, according to conjecture \ref{c}, the logarithmic Kobayashi pseudometrics of these subvarieties should be generically metrics. A statement in this direction is shown in \cite{VZ01}. If one could prove that any generic two points of an orbifold can be connected by chains of orbifold entire curves $h:\Bbb C\to (X\vert \Delta)$ (see remark \ref{cs}.(5)), the conjecture would also follow from \cite{VZ01}.

3. Since the conjectures  \ref{abcj},\ref{lmmp} and \ref{c1=0} are known in dimension at most $2$, and also in dimension $3$ if $\Delta=0$, we see that the isotriviality conjecture is true when $d=3$ (as shown in \cite{JK}), and also when $d=4$ if $\Delta=0$. 
\end{remark}

\section{\bf Special versus Weakly Special Manifolds}\label{wspec}

In this section, we illustrate by an example the difference between the two notions, and its implications in arithmetics and hyperbolicity questions. It also shows that orbifold structures need to be considered in birational classification.

Recall (definition \ref{ws}) that the (complex, projective, connected) manifold $X$ is `weakly special' if none of its finite \'etale covers maps rationally onto a positive-dimensional variety of general type. Special $X's$ are weakly special, and conversely for curves and surfaces.

However simply-connected threefolds which are weakly special, but not special, are constructed in \cite{BT} (see remark \ref{bt} above). Some of their examples are defined over $\Q$. We have the following conjecture, due to Abramovich and Colliot-Th\'el\`ene, stated in \cite{HT}:

\begin{conjecture}\label{act} (Abramovich-Colliot-Th\'el\`ene) Let $X$ be defined over a number field. Then $X$ is potentially dense if and only if $X$ is weakly special (as formulated in our terminology).
\end{conjecture}

The `only if' part follows from Lang's conjecture and Chevalley-Weil theorem.

Observe that this conjecture conflicts with item $3$ of conjecture \ref{cspec} for Bogomolov-Tschinkel threefolds $X$ (if defined over $\Q$, say). Indeed, \ref{act} above claims that $X$ is potentially dense, while our conjecture \ref{cgen} claims they are not (and even that $c_X(X(k))\cap U$) is finite for any number field $k$.

 Let us sketch the construction of the Bogomolov-Tschinkel threefolds given in \cite{BT}.

\begin{example}\label{btex} One needs the following ingredients:

1. A smooth elliptic surface $g:T\to\PP^1$ with exactly one multiple fibre $m.T_0=g^{*}(0)$ of multiplicity $m\geq 2$, with $T_0$ smooth connected, such that $(T-T_0)$ is simply-connected. Such a surface can be constructed by a logarithmic transform on a suitable elliptic rational surface.

2. A smooth simply-connected elliptic surface $h':S'\to \PP^1$ with $\kappa(S')=1$, together with an ample base-point free line bundle $L'$, and a generic pencil of (generically smooth) sections of $L'$ giving a map $\varphi: S'\dasharrow \PP^1$. By Lefschetz theorem, the complement in $S'$ of a smooth member of this pencil has a finite fundamental group, equal to the order of divisibility of $L'$ in $Pic(S')$.

Two generic members $B,B'$ of the pencil meet transversally at a finite set $P'$ of points $p'_j$ of $S$. Blow-up these points to get a regular map $\psi: S\to \PP^1$ on the blown-up surface $S$. Denote by $D,D'$ the strict transforms of $B,B'$ in $S$.

The crucial observation is that their complements in $S$ are simply-connected (this is because the lift to $S$ of a small loop around $B$ in $S'$, and close to $p_j$ say, becomes homotopically trivial in $S$, since already homotopically trivial in $(E_j-q_j)$, if $E_j$ is the exceptional divisor of $S$ above $p_j$, and $q_j$ is the intersection point of $D$ with $E_j$)

The sought-for elliptic threefold $f:X\to S$ is then just  obtained from $\psi:S\to \PP^1$ by the base change $g:T\to \PP^1$. It has a multiple fibre of multiplicity $m$ exactly above the fibre $D$ of $\psi$, so that $\Delta(f)=(1-\frac{1}{m}).D$. From which one immediately deduces that $\kappa(S)=1<\kappa(S,K_S+\Delta_f)=2$. It only remains only to check that $X$ is simply connected, which easily follows from the simple-connecteness of $(T-T_0)$ and $(S-D)$ $\bullet$

\end{example}

The known methods of arithmetic geometry are presently unable to solve the problem of whether these threefolds are potentially dense or not (One had to decide whether certain simply connected smooth orbifold surfaces of general type are potentially dense or not). Instead of this, one can answer the hyperbolic analogue for some examples at least. Let us state the hyperbolic analogue of Abramovich-Colliot-Th\'el\`ene conjecture \ref{act}:

\begin{question} Is a (complex, projective, connected) manifold $X$ `weakly special' if and only if $d_X\equiv 0$?
\end{question}

The answer to this question is `no' for the `only if' part (the `if' part follows from Lang's conjecture). This shows that either conjecture \ref{act} is wrong, or that the expected links between arithmetics and hyperbolicity fail to hold. More precisely:

\begin{theorem}(\cite{CP}) There exists certain simply connected elliptic threefolds $c_X:X\to S$ constructed in \cite{BT} (and so: weakly special, but not special) such that, for any entire curve $h:\C\to X$, the composed map $c_X\circ h:\C\to S$ is either constant, or has image contained in some fixed projective curve $C\subset S$ (independent on $h)$.

(\cite{CW}) Certain of the above examples can be chosen so that the curve $C$ is empty. In this case, $d_X=(c_X)^*(\delta)$, for some continuous metric $\delta$ on $S$ (which establishes conjecture \ref{cgen} in these cases).
\end{theorem}

The proof consists in adapting to the orbifold context the method of Bogomolov to show that surfaces of general type with $(c_1^2-c_2)>0$ contain only finitely many special curves, as extended to the transcendental case by Mc Quillan. The methods are applied to the (general type) orbifold base $(S\vert \Delta_{c_X})$ of $c_X$ to show that all orbifold maps $h':\C\to (S\vert \Delta_{c_X})$ are algebraically degenerate. Because $f\circ h$ is such an orbifold map for any $h:\C\to X$, the conclusion follows $\square$

\

In \cite{Rou06} this result is generalised, and its proof simplified, using some further results of Mc Quillan.

\section{\bf Classical versus `non-classical' multiplicities}\label{classmult}

The introduction of multiplicities defined by infimum rather than gcd creates great technical difficulties, although possibly not modifying the qualitative geometry. We illustrate this by two examples from \cite{Ca05}. See also conjecture \ref{zvsqrat}. 

We first start (see \ref{tncf} below) with an example showing a huge discrepancy in the case of fibres of general type case. We ask then (see questions \ref{qcvsnc} and \ref{qcvsnc'}) whether this discrepancy can occur when the fibres are special.

\begin{theorem}\label{tncf} There exists {\bf general type} fibrations $f:S\to \PP^1$ with $S$ a smooth projective connected and {\bf simply-connected} surface.
\end{theorem}

Note that the orbifold base cannot be the `classical' (or `divisible') orbifold base. Indeed, if $f:S\to \PP^1$ were of general type for the `classical' multiplicities, a finite \'etale cover of $S$ would map onto a hyperbolic curve, and its fundamental group had the free group on two generators as a quotient. 
Note also that the smooth fibres of $f$ are hyperbolic curves, since multiple fibres of elliptic surfaces are always `divisible' (ie: $inf$ and $gcd$ coincide for them). In fact in the examples of \cite{Ca05}, the genus f the fibres is $23$. Other more interesting examples with fibres of any genus at least $2$ have been constructed by Lidia Stoppino (\cite{Sto}, unfortunately unpublished, but going to be put arXiv). Finally, observe that $S$ is necessarily of general type (by orbifold additivity, for example).

$\bullet$ Let us give a brief sketch of the construction of \cite{Ca05}: for suitable choice of five distinct  lines $T_k,k=1,2$ and $D_j,j=1,2,3$ of $\PP^2$ meeting in a point $a$, we assign to the lines $D_j$ the multiplicity $2$, and to the lines $T_k$ the multiplicity $3$ getting a (non-reduced) curve $C$ of degree $2.3+3.2=12$, we show the existence of an irreducible  curve $C'\subset \PP^2$ of degree $12$ not going through $a$, and meeting each of the $T_k$ (resp. $D_j$) in $4$ (resp. $6$) distinct points with order of contact exactly $3$ (resp. $2$). The pencil of curves of degree $12$ on $\PP^2$ generated by $C,C'$ gives a rational map $h':\PP^2\to \PP^1$ which we resolve as $h:P\to \PP^1$. One checks then that it has a (non-classical) multiple fibre $C_0$ (of multiplicity $2$) consisting of a (simply connected) tree of rational curves with multiplicities either $2$ or $3$. To get $f: S\to \PP^1$, one just make base change through a map: $r:\PP^1\to \PP^1$ of degree at least $5$ ramified at generic points. The simple-connectedness of $S$ follows from the exact sequence of groups: $\pi_1(F)\to \pi_1(S)\to \pi_1(\PP^1)$, if $F$ is a generic fibre of $f$. The exactness of this sequence is because $f$ has no classical multiple fibre, the image of $\pi_1(F)$ in $\pi_1(S)$ is trivial because $f$ has some simply connected fibres (like $C_0$) $\bullet$

\

Some of the examples above are defined over $\bQ$. For any number field $k$, $f(S(k))$ should be finite, after conjecture \ref{cgen}, which would establish Lang's conjectural `mordellicity'\footnote{ie: $S(k)\cap U$ is finite for any $k$ and some fixed Zariski open nonempty $U\subset S$.} for them, and give the first example of a `mordellic' smooth simply connected surface of general type. Because $f(S(k))$ is contained in the set of rational points of the base orbifold $(\PP^1\vert \Delta_f)$ of $f$, which is of general type, the problem were solved if one could establish the orbifold Mordell conjecture \ref{com} below.

\begin{question}\label{qcvsnc} Let $f:X\to B$ be a fibration with $X$ smooth, and generic (smooth) fibres $F$ having $\kappa(F)=0$. Do we have the equality: $\Delta_f=\Delta^*_f$? (In other words: do the $inf$ and $gcd$-multiplicities then coincide ?).

This holds when $F$ are abelian varieties. The first non-trivial case is for $K3$ surfaces (and even Kummer surfaces).
\end{question}

More generally:

\begin{question}\label{qcvsnc'}

1. Let $f:X\to B$ be a fibration with $X$ smooth, and generic (smooth) fibres special. Do we have the equality: $\Delta_f=\Delta^*_f$?
(Observe indeed that the property holds  when $F$ is rationally connected, by \cite{GHS}).

2. If $X$ is a `classical' special manifold, is it special? (Being `classically' special means that there is no neat (rational) fibration on $X$ with classical orbifold base of general type).

\end{question}

\section{\bf An orbifold version of Mordell's conjecture}

 \
 
We shall introduce two notions (`classical' and `non-classical') of $\Q$-rational points on certain orbifolds of the form $(\PP^1\vert \Delta)$. These two notions are deduced from the two natural notions (based either on $gcd$ or $inf)$ of orbifold morphisms, when rational points are seen as sections of the arithmetic surface $\PP^1$ over the spectrum of the ring of integers.

These two notions are also respectively compatible functorially with the two natural notions of fibre multiplicities considered in the text in the following sense: if $f:S\to \PP^1$ is a fibration with orbifold base (resp. `classical' orbifold base) $(\PP^1\vert \Delta)$, everything defined over $\Q$, then $f(X(\Q))$ is contained in the set of $\Q$-rational (resp. `classical' $\Q$-rational) points of $(\PP^1\vert \Delta)$ defined below.

 \
 
  Let $(\PP^1\vert\Delta)$ be a geometric orbifold, with integers $p,q,r>1$, and: 
  
  $\Delta:=(1-1/p).\{0\}+(1-1/q).\{1\}+(1-1/r).\{\infty\}$.
  
  Then: $\kappa(\PP^1\vert\Delta)=1$  if and only if: $(1/p+1/q+1/r)<1$. 
  
  This is the case for example if: $(p,q,r)=(2,3,7)$
  
  This geometric orbifold is defined over $k=\Bbb Q$. 
  
  Let us now define its $\Q$-rational points, and also their `classical' version.
  
  \

  $\bullet$ The set  $(\PP^1\vert\Delta)(\Q)^*$ of `classical' $\Q$-rational points of $(X\vert\Delta)$ consists of the usual rational points $x=\frac{\alpha^p}{\beta^r}$, with $\alpha,\beta\in \Z, gcd(\alpha,\beta)=1$, such that there exists: $\gamma\in \Z$, and: $ \alpha^p+\beta^r=\gamma^q$.

  These points are the sections of the arithmetic surface $\PP^1$ over the spectrum of the ring of integers, meeting the three sections defined by the points $0,\infty,1$ with arithmetic orders of contact divisible by $p$, $q$ and $r$ respectively.

  Using Falting's and Chevalley-Weil's theorems, one gets:

\begin{theorem}(\cite{DG})( `Faltings+$\varepsilon$') 

$(\PP^1\vert\Delta)(\Q)^*$ is finite  if $\kappa(\PP^1\vert\Delta)=1$.
  \end{theorem}

  $\bullet$ The set $(\PP^1\vert\Delta)(\Q)$ of (`non-classical') $\Q$-rational points of $(\PP^1\vert\Delta)$ are the usual rational points $x=\frac{a}{b}$ with $a,b\in \Z, (a,b)=1$, such that:

  $a$ is $`p$-full',$b$ is $`r$-full', and $c:=a-b$ is $`q$-full', where:
  
   $a$ is $`p$-full' means: for any prime $\ell$ dividing $a$, then: $\ell^p$ divides $a$. 
   
  Notice that, by a result of Erd\" os, the set of $p$-full integers behaves asymptotically as the set of $p$-th powers: $Card\{a\leq X\vert a$ is $p$-full$\}\sim_{X\to+\infty}C_p.X^{1/p}$, for a certain explicit constant $C_p>0$.
  
   These points are the sections of the arithmetic surface $\PP^1$ over the spectrum of the ring of integers, meeting the three sections defined by the points $0,\infty,1$ with arithmetic orders of contact {\it at least} equal to $p$, $q$ and $r$ respectively.

  \begin{conjecture}\label{com} (Orbifold Mordell conjecture)

  $(\PP^1\vert\Delta)(\Q)$ is finite if $\kappa(\PP^1\vert\Delta)=1$
  \end{conjecture}

   \begin{remark} The above conjecture is open. It is easy to show\footnote{This observation, due to P. Colmez, has been communicated to me by Colliot-Th\'el\`ene.} that the $abc$-conjecture implies the  Orbifold Mordell conjecture. This conjecture can be stated similarly (but less concretely) for arbitrary number fields $k$. It is a (very) special case of the conjecture \ref{cgen} above.
   \end{remark}

\end{document}